\documentclass{amsart}

\usepackage[utf8]{inputenc}

\usepackage{mathtools}
\usepackage{amssymb}
\usepackage{amsthm}
\usepackage{amsmath}
\usepackage{bm}
\usepackage{enumitem}
\usepackage{cite}
\usepackage{xfrac}
\usepackage{xcolor}
\usepackage{tikz-cd}
\usepackage{lipsum}
\usepackage{adjustbox}
\usepackage{pdfpages}
\usepackage{hyperref}

\theoremstyle{plain}
\newtheorem{theorem}{Theorem}[section]
\newtheorem*{theorem*}{Theorem}
\newtheorem{corollary}[theorem]{Corollary}
\newtheorem*{corollary*}{Corollary}
\newtheorem{lemma}[theorem]{Lemma}
\newtheorem*{theoremA*}{Theorem A}
\newtheorem*{theoremB*}{Theorem B}
\newtheorem*{theoremC*}{Theorem C}

\theoremstyle{definition}
\newtheorem{definition}[theorem]{Definition}
\newtheorem{examples}[theorem]{Examples}
\newtheorem*{definition*}{Definition}

\newtheorem{fact}[theorem]{Fact}

\theoremstyle{remark}
\newtheorem{remark}[theorem]{Remark}

\theoremstyle{remark}
\newtheorem{notation}[theorem]{Notation}
\newtheorem*{claim}{Claim}

\theoremstyle{plain}
\newtheorem{proposition}[theorem]{Proposition}

\theoremstyle{plain}
\newtheorem{diagram}[theorem]{Diagram}

\newcommand{\VF}{\mathrm{VF}}
\newcommand{\VG}{\mathrm{VG}}
\newcommand{\RF}{\mathrm{RF}}

\newcommand{\rplus}{(\mathbb{R}^{+},\cdot)}

\newcommand{\X}{\Bar{X}}
\newcommand{\Xs}{\Bar{X}^{*}}

\newcommand{\PX}{P(\Bar{X})}
\newcommand{\Phx}{P^{h}(\Bar{X},\Bar{X}^{*})}
\newcommand{\Pxx}{P(\Bar{X},\Bar{X}^{*})}
\newcommand{\Ps}{P^{*}(\Bar{X}^{*})}
\newcommand{\PSxx}{P_{S}(\Bar{X},\Bar{X}^{*})}

\newcommand{\Zx}{\mathbb{Z}[\Bar{X}]}
\newcommand{\Zhx}{\mathbb{Z}^{h}[\Bar{X}]}
\newcommand{\Zxx}{\mathbb{Z}[\Bar{X},\Bar{X}^{*}]}
\newcommand{\Zxn}{\mathbb{Z}[X_{1},\dots,X_{n}]}

\newcommand{\preds}[1]{\|#1^{*}\|}

\newcommand{\predpnxb}{\|P_{n}(\Bar{x})\|}
\newcommand{\predpsxb}{\|P^{*}(\Bar{x})\|}

\newcommand{\predpab}{\|P(\Bar{\am})\|}

\newcommand{\predbs}{\|\bmodel^{*}\|}

\newcommand{\predxs}{\|x^{*}\|}
\newcommand{\predxis}{\|x_{i}^{*}\|}

\newcommand{\Ov}{\mathcal{O}_{v}}
\newcommand{\Ow}{\mathcal{O}_{w}}

\newcommand{\valuedfield}[1]{(#1,\Gamma_{#1},k_{#1})}
\newcommand{\Kfield}{(K,\Gamma_{K},k_{K})}
\newcommand{\mv}{\textit{m}_{v}}
\newcommand{\pl}{\text{pl}}

\newcommand{\af}{a}
\newcommand{\bfield}{b}
\newcommand{\cf}{c}

\newcommand{\kpro}{K\mathbb{P}^{1}}
\newcommand{\fpro}{F\mathbb{P}^{1}}

\newcommand{\upkipmetr}{(\prod_{i\in\mathcal{I}}K_{i}\mathbb{P}^{1})_{\mathcal{D}}^{\mathrm{me}}}

\newcommand{\upkicl}{(\prod_{i\in\mathcal{I}}k_{i})_{\mathcal{D}}^{\mathrm{cl}}}
\newcommand{\upgicl}{(\prod_{i\in\mathcal{I}}\Gamma_{i})_{\mathcal{D}}^{\mathrm{cl}}}
\newcommand{\upKicl}{(\prod_{i\in\mathcal{I}}K_{i})_{\mathcal{D}}^{\mathrm{cl}}}
\newcommand{\upKcalcl}{(\prod_{i\in\mathcal{I}}\mathcal{K}_{i})_{\mathcal{D}}^{\mathrm{cl}}}

\newcommand{\prodki}{\prod_{i\in\mathcal{I}}K_{i}}
\newcommand{\prodkip}{\prod_{i\in\mathcal{I}}K_{i}\mathbb{P}^{1}}

\newcommand{\limiD}{\lim_{i\rightarrow\mathcal{D}}}

\newcommand{\Rfin}{\mathbb{R}^+_{\mathrm{fin}}}
\newcommand{\gfin}{\Gamma_{\mathrm{fin}}}
\newcommand{\Rinf}{\mathbb{R}^+_{\mathrm{inf}}}
\newcommand{\ginf}{\Gamma_{\mathrm{inf}}}
\newcommand{\gfininf}{\sfrac{\gfin}{\ginf}}

\newcommand{\vinf}{v_{\mathrm{inf}}}
\newcommand{\vfin}{v_{\mathrm{fin}}}
\newcommand{\vfinb}{\bar{v}_{\mathrm{fin}}}

\DeclareMathOperator{\plvb}{pl_{\bar{\mathrm{v}}}}

\DeclareMathOperator{\dg}{\mathrm{dg}}
\newcommand{\dgk}{\dg_{K}}
\newcommand{\dgf}{\dg_{F}}

\newcommand{\am}{\mathbf{a}}
\newcommand{\bmodel}{\mathbf{b}}

\newcommand{\first}[1]{{#1}^{\circ}}
\newcommand{\ao}{a^{\circ}}
\newcommand{\bo}{b^{\circ}}

\newcommand{\second}[1]{{#1}^{*}}
\newcommand{\as}{a^{*}}
\newcommand{\bs}{b^{*}}

\newcommand{\n}[1]{\|#1\|}
\newcommand{\sn}[1]{\|{#1}^{*}\|}

\newcommand{\mvf}{\mathrm{MVF}_{0,0}}
\newcommand{\mvfd}{\mathrm{MVF}_{0,0}^{\mathrm{d}}}

\newcommand{\lp}{\mathcal{L}_{\mathbb{P}}}
\newcommand{\lps}{\mathcal{L}_{\mathbb{P},\sigma}}
\newcommand{\lval}{\mathcal{L}_{\mathrm{val}}}
\newcommand{\lvals}{\mathcal{L}_{\mathrm{val},\sigma}}

\newcommand{\lcval}{\mathcal{L}_{\mathrm{c,val}}}

\newcommand{\lring}{\mathcal{L}_{\mathrm{ring}}}
\newcommand{\lordgr}{\mathcal{L}_{\mathrm{og}}}
\newcommand{\tcval}{T_{\mathrm{c,val}}}

\def\dotminussym#1{%
  \setbox0=\hbox{$-$}%
  \kern.5\wd0%
  \hbox to 0pt{\hss\hbox{$-$}\hss}%
  \raise.6\ht0\hbox to 0pt{\hss$.$\hss}%
  \kern.5\wd0%
}

\title{An Approximate AKE Principle for Metric Valued Fields}

\author{Martin Hils}
\address{Institut f\"{u}r Mathematische Logik und Grundlagenforschung, Universit\"{a}t M\"{u}nster, Einsteinstr. 62, D-48149 M\"{u}nster, Germany}
\email{hils@uni-muenster.de}
\thanks{MH was partially supported by the German Research Foundation (DFG) via HI 2004/1-1 (part of the French-German ANR-DFG project GeoMod) and under Germany's Excellence Strategy EXC 2044-390685587, `Mathematics M\"unster: Dynamics-Geometry-Structure'.}

\author{Stefan Marian Ludwig}
\address{Albert-Ludwigs-Universität Freiburg,
Mathematisches Institut,
Abteilung für Mathematische Logik,
Ernst-Zermelo-Straße 1,
79104 Freiburg i. B., Germany}
\email{stefan.ludwig@mathematik.uni-freiburg.de, ludwig@dma.ens.fr}	
\thanks{SML has received funding from the European Union's Horizon 2020 research and innovation programme under the Marie Sk\l{}odowska-Curie grant agreement N\textsuperscript{\underline{o}} 945322.  \includegraphics[scale=0.025]{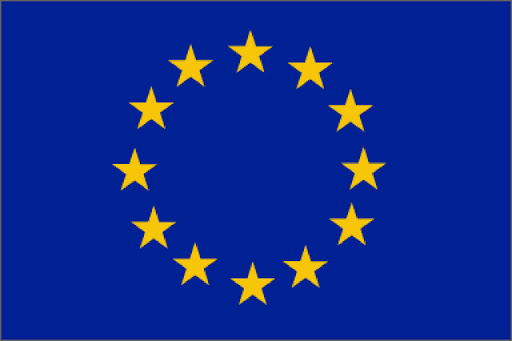}	Moreover, he was partially supported by GeoMod AAPG2019 (ANR-DFG), `Geometric and Combinatorial Configurations in Model Theory', and as a 2-month stagiaire at the Institut Camille Jordan (UCBL)}



\begin{document}

\begin{abstract} 
We study metric valued fields in continuous logic, following Ben Yaacov's approach, thus working in the metric space given by the projective line. As our main result, we obtain an approximate Ax-Kochen-Ershov principle in this framework, completely describing elementary equivalence in equicharacteristic 0 in terms of the residue field and value group. Moreover, we show that, in any characteristic, the theory of metric valued difference fields does not admit a model-companion. This answers a question of Ben Yaacov.
\end{abstract}

\maketitle


\section{Introduction}
In \cite{YAACOV_2014} Ben Yaacov introduced a formalism to consider certain valued fields, called \textit{metric valued fields}, as structures in continuous logic, namely complete valued fields with value group embedded in the group $\rplus$. Given a metric valued field $K$, for technical reasons Ben Yaacov associates to it a continuous logic structure $K\mathbb{P}^1$ with base set the projective line over $K$.\footnote{Very recently, Rideau-Kikuchi, Scanlon and Simon used an alternative approach in \cite{RiSc25}.}
Still in \cite{YAACOV_2014},  he further established a quantifier elimination result for the theories of (projective lines over) algebraically closed and of  real closed metric non-trivially valued fields.

In this  article,  we will deepen the study of metric valued fields in this context. Naturally the question arises, whether there exists a connection between residue field and value group of a metric valued field in  equicharacteristic  $0$ and its elementary (continuous logic) theory, similarly to the Ax-Kochen-Ershov (AKE) principle in the classical context. This question is  non-trivial  since in general the residue field of $K$ is not interpretable in the continuous logic structure $K\mathbb{P}^1$ and the value group is directly given by the metric.  Consequently, we not only have to investigate to which extent the elementary theories of residue field and value group determine the elementary theory of the metric valued field but also vice-versa, unlike in the classical context.  In our analysis, we  restrict  ourselves mainly to the case of equicharacteristic $0$ metric valued fields with \textit{dense} (in $\rplus$) value group, as non-dense metric valued fields resemble classical logic structures.

As our main result, we obtain an \textit{approximate AKE principle}, showing that the elementary theory of a metric valued field of equicharacteristic $0$ with dense value group determines and is determined by the elementary theories of its residue field  and value group, up to what we will call a  \textit{residue shift}. For this,  we introduce  the classes $\mathcal{C}(\mathfrak{G},\mathfrak{F})$. The idea is that for countably incomplete metric ultrapowers of metric valued fields as above the residue field carries an infinitesimal valuation and the potentially changing elementary theory of the residue field is still controllable since this infinitesimal valuation is determined up to elementary equivalence. 
We call a pair $(\mathfrak{G},\mathfrak{F})$ admissible, if $\mathfrak{G}$ is the $\mathcal{L}_{\mathrm{og}}$-theory of a regular dense (non-trivial) ordered abelian group $\Delta$ and $\mathfrak{F}$ is the $\lring$-theory of a field $l$ of characteristic $0$ and such that the pair $(\Delta, l)$ belongs to one of the following cases: \begin{itemize}
    \item[(i)] $\Delta$ is not divisible; 
    \item[(ii)] $\Delta$ is divisible and $l\not\equiv l^{\prime}((t^{\Delta^{\prime}}))$ in $\lring$ for any such pairs $(\Delta^{\prime},l^{\prime})$ with $\Delta^{\prime}$ non-divisible or with $\Delta'$ divisible and $l\not\equiv l^{\prime}$ (see Definition~\ref{def:generating-pair}).
\end{itemize}
Given an admissible pair $(\mathfrak{G},\mathfrak{F})$ we define the class $\mathcal{C}(\mathfrak{G},\mathfrak{F})$ to consist of all metric valued fields $\mathcal{K}=(K,\Gamma_{K},k_{K})$ 
such that $\Gamma_{K}\models \mathfrak{G}$ and $k_{K}\models \mathfrak{F}$ or $\Gamma_{K}\models \mathrm{DOAG}$ and $k_{K}\models \mathrm{Th}(l((t^{\Delta})))$. Now, the main result reads as follows.
\begin{theoremA*}[Theorem~\ref{metricAKEnew}]
Let $K,F$ be metric valued fields of equicha\-rac\-teristic $0$
 with dense value groups. Then $\kpro$ and $\fpro$ are elementarily equivalent if and only if they are in the same class $\mathcal{C}(\mathfrak{G},\mathfrak{F})$ for some admissible pair $(\mathfrak{G},\mathfrak{F})$.
\end{theoremA*}
To establish this result, a key tool will be to relate the metric ultrapowers of metric valued fields to their classical logic counterparts.  We will build on the fact, that the metric value group is always archimedean  and hence  regular   in the sense of \cite{Elt.-prop.-of-ordered-abelian-groups}. Moreover, due to a result of Hong \cite{defbilityHong} (building on Koenigsmann \cite{koenigsmann}), the valuation is then definable (in classical logic) in the non-divisible case. 
As a consequence of our main theorem,  we can retrieve the following as a corollary, showing that despite the existence of the valuation and  its  presence  as  the underlying metric,  elementary equivalence in the continuous context in the end reduces to elementary equivalence in the language $\lring$.

\begin{theoremB*}[Corollary~\ref{Cor:1storder}]
Let $K,F$ be metric valued fields of equicha\-rac\-teristic $0$
 with dense value groups. Then $\kpro\equiv\fpro$ if and only if $K\equiv F$ in $\lring$.
 \end{theoremB*}

The model theory of fields enriched with an automorphism proved to be an interesting object of study leading to striking applications, e.g., in diophantine geometry \cite{HRUSHOVSKImaninmumford} and in algebraic dynamics \cite{chatzidakis_hrushovski_2008,medvedevscanlon}. The theory of difference fields (fields endowed with an automorphism) admits a model-companion $\mathrm{ACFA}$, a simple unstable theory which has been developed in fundamental work by Chatzidakis-Hrushovski \cite{Difference} and which lies at the heart of these applications.

By a classical result of A. Robinson, the theory $\mathrm{ACVF}$ of algebraically closed non-trivially valued fields is the model-companion of the theory of valued fields in classical logic. Ben Yaacov obtained the analogous result in the metric setting. While all completions of $\mathrm{ACVF}$ are NIP and unstable, the model-companion $\mathrm{ACMVF}$ of the theory of metric valued fields is stable.

By \cite{Kikyo2002-KIKTSO} the theory of valued fields with an automorphism does not admit a model-companion, even when restricted to equicharacteristic 0. Although, by requiring the induced automorphism on the value group to be  the identity,  one may overcome this.  This theory does admit a model-companion in equicharacteristic 0, by work of B\'elair-Macintyre-Scanlon \cite{frobeniuswitt}, whose completions are $\mathrm{NTP}_2$, as shown in \cite{ChHi14}.

Ben Yaacov asked in 2018 whether the theory of (projective lines over) metric valued fields in equicharacteristic 0, endowed with an isometric automorphism, admits a model-companion, which is a very natural question given the results we mentioned.  If the model-companion existed, it would be a candidate for a natural continuous-logic example of a simple unstable theory. In our paper, we answer this question. Rather surprisingly, we show that a model-companion does not exist in the continuous context. Concretely, we obtain the following result.

\begin{theoremC*}[Theorem~\ref{thm:no-model-companion}]
Fix any $(a,b)\in\{(0,0),(0,p),(p,p)\mid p\textrm{ prime}\}$. Then the theory of \textit{metric valued difference fields} of characteristic $(a,b)$  does not have a model-companion.
\end{theoremC*}
The key idea behind the proof is, that whereas in the context without automorphism the phenomenon of the residue shift is  controllable,  the interaction of residue shift and  automorphism  is not.

\subsection*{Structure of the article} 
In Section~\ref{S:prelim} we first recall some facts from the classical model theory of valued fields and then focus on regular (ordered) abelian groups.
Here, we present Hong's definability result  
and deduce the consequences we will use in the proof of the main theorem. Further, we give an introduction to the continuous logic used and to Ben Yaacov's formalism 
to handle metric valued fields in continuous logic. In Section~\ref{S:ult} we investigate ultraproducts of metric valued fields and obtain what we call the \textit{residue shift}. In Section~\ref{S:main} we state the main theorem together with some examples, while Section~\ref{S:proof} is devoted to its proof. Finally, in Section~\ref{S:diff}, we prove Theorem~C.

\subsection*{Acknowledgements} The results of this article were obtained in the master thesis of the second author at the University of M\"unster, supervised by the first author. The second author would like to thank Ita\"{i} Ben Yaacov for introducing him to the topic and many helpful discussions during an extended research stay at Universit\'e Claude Bernard Lyon 1.

\section{Preliminaries}\label{S:prelim}
\subsection{Notation and Prerequisites}
We will briefly recall some notations and results from the model theory of valued fields. For an in-depth treatment of valued fields we refer the reader to \cite{engler2005valued}. As we will later work in a metric context the valuations will be written mutliplicatively throughout the paper. In the following by an ordered group we will always mean a totally ordered group. 
\begin{definition}
Let $K$ be a field. Then $(K,v)$ is called a \textit{valued field} if $v:K\rightarrow \Gamma \cup \{0 \}$ is a \textit{valuation map} (or simply \textit{valuation}), i.e., we have that $v(x)=0$ if and only if $x=0$ and further $v(xy)=v(x)\cdot v(y)$ and $v(x+y)\leq \max\{v(x),v(y)\}$ holds for all $x,y\in K$. Here, 
$\Gamma=(\Gamma,\cdot,1,<)$ is required to be an ordered abelian group with $0<\Gamma$ and it will be called the \textit{value group}.
Further, $\Ov:=\{x\in K\;|\;v(x)\leq 1\}$ is the \textit{valuation ring} of $K$ with maximal ideal $\mv:=\{x\in K\;|\;v(x)<1\}$.

The quotient $\sfrac{\Ov}{\mv}$ is called the \textit{residue field} of $K$, mostly denoted by $k$, $k_{v}$ or $k_{K}$. We write $\text{res}:\Ov\rightarrow k$ for the quotient map which is called the $\textit{residue map}$.
\end{definition}
\begin{definition}
Given fields $K$ and $k$, a \textit{place} $\pl:K\rightarrow k\cup\{\infty\}$ is a surjective map such that the following properties hold:
\begin{itemize}
    \item If $x\in K^\times$ then $\pl(x)=0$ if and only if $\pl(1/x)=\infty$.
    \item If $x,y\in K$ then $\pl(x+y)=\pl(x)+\pl(y)$ and 
    $\pl(xy)=\pl(x)\pl(y)$ hold whenever they are defined using the conventions that $\infty+a=\infty=a+\infty$ for all $a\in k$ and $\infty\cdot a=\infty=\infty\cdot\infty=a\cdot\infty$ for all $a\in k^{\times}$. The expressions $\infty\cdot 0$, $0\cdot\infty$ and $\infty+\infty$  are not defined.
\end{itemize}
\end{definition}

\begin{remark}
Recall that there is a 1:1-correspondence between valuation rings of a field $K$ and valuations of $K$ up to equivalence. Concretely, a valuation $v$ on $K$ gives rise to the valuation ring $\Ov$, and a valuation ring $\mathcal{O}$ determines a valuation $v$ with values in the abelian group $K^{\times}/\mathcal{O}^{\times}$ ordered by setting $x\mathcal{O}^{\times}\leq y\mathcal{O}^{\times}$ whenever $xy^{-1}\in\mathcal{O}$, where $v$ is defined on $K^{\times}$ by $v(x):=x\mathcal{O}^{\times}$ and $v(0):=0$. Moreover we have the following connection between places and valuations:
A valuation $v:K\rightarrow\Gamma$ with residue map $\text{res}:\Ov\rightarrow k$ determines a place $\pl_{v}:K\rightarrow k\cup\{\infty\}$ by setting $\pl_{v}(x)= \text{res}(x)$ if $x\in\Ov$ and $\pl_{v}(x)=\infty$ otherwise.
On the other hand given a place $\pl:K\rightarrow k$ we can deduce a corresponding valuation by setting $\pl^{-1}(k)$ to be the valuation ring with maximal ideal $\pl^{-1}(\{0\})$.
\end{remark}

Recall that a valued field $\Kfield$ is of \textit{equicharacteristic 0} (resp.\ $p$) if $\text{char}(K)=0=\text{char}(k_{K})$ (resp.\ $\text{char}(K)=p=\text{char}(k_{K})$). If $\text{char}(K)=0$ but $\text{char}(k_{K})=p$ for some prime $p$ one says $K$ is of \textit{mixed characteristic (0,p)}.

\begin{definition}
The language $\lval$ consists of a sort $\VF$ for the main field in the language of rings $\lring$, a sort $\VG$ for the value group in $\lordgr\cup\{0\}$, the language of ordered groups with a constant symbol for $0$ and a sort $\RF$ for the residue field in $\lring$ again. Furthermore $\lval$ contains function symbols $v:\VF\rightarrow \VG$ and $\mathrm{Res}: \VF^{2}\rightarrow \RF$.\\
In a valued field $v$ shall be interpreted as the valuation and
\[\mathop{\mathrm{Res}}(x,y):=\begin{cases}
\text{res}(xy^{-1}), & \text{ if }0<v(x)\leq v(y), \\
0, & otherwise.
\end{cases}\]
The theory of henselian valued fields of equicharacteristic $0$ (with the above interpretations) in $\lval$ shall be denoted by $T_{\mathrm{val}}$. \\
Moreover the language $\lcval$ shall consist of an additional constant symbol $c$ in the value group sort and $\tcval$ shall contain an additional axiom stating that neither $c$ is the identity element of the value group, nor $c=0$.
\end{definition}
Now we state the so-called \textit{AKE principle} which will play a significant role throughout this paper. It was first obtained by Ax and Kochen \cite{diohantineProblems, AK(E)2,AK(E)3} and Er\v{s}ov \cite{(AK)E} independently. Note that the statement for $\lcval$ can be easlily obtained from the statement in $\lval$ using the pure stable embeddedness of the value group.
\begin{theorem}\label{AKE}
Let $\mathcal{K}=\Kfield$ and $\mathcal{F}=\valuedfield{F}$ be henselian valued fields in equicharacteristic $0$. Then the following holds in $\lval$ (resp. $\lcval$):
\[\mathcal{K}\equiv \mathcal{F} \;\;\text{if and only if}\;\; k_{K}\equiv k_{F} \;\text{and}\;\Gamma_{K}\equiv\Gamma_{F}.
\]
\end{theorem}

\subsection{Regular valuations}
Value groups of metric valued fields will be subgroups of $\mathbb{R}^{+}=(\mathbb{R}^{+},\cdot,1,<)$ and consequently archimedean ordered abelian groups. While archimedeanity itself is not axiomatizable in classical discrete model theory one can consider a larger class of groups called regular ordered abelian groups. This class  was introduced by A. Robinson and Zakon in \cite{Elt.-prop.-of-ordered-abelian-groups} and further elaborated in \cite{generalized_archimedean_groups}. It is indeed an elementary class.\\
Returning to fields, it is in general an interesting question to ask when the existence of a (henselian) valuation already implies the definability of this valuation in the language $\lring$, i.e., that the valuation ring is a definable set. When determining the complete theory of a metric valued field in terms of its residue field and value group we will later often deal with the case of a residue field that itself carries a valuation with regular value group. In order to capture this valuation we will make use of a definability result 
as given by Hong \cite{defbilityHong}.\\
As Definition \ref{definitionregulargroup} to Theorem \ref{convex Sg. elt. equiv.} only deal with the ordered (abelian) group structure we use, different from the rest of the paper, an additive notation. 

\begin{definition}\label{definitionregulargroup}
A non-trivial ordered abelian group $(G,+,<)$ is called 
\begin{itemize}
    \item \textit{discrete} if it has a minimal positive element, and \textit{dense} otherwise;
    \item \textit{regular}, if for any non-trivial convex subgroup $H\subseteq G$ the quotient group $G/H$ is divisible;
    \end{itemize}
    A valuation with regular value group is called a \textit{regular} \textit{valuation}.
\end{definition}

\begin{remark}
    The theory of ordered abelian group has a model completion, denoted by DOAG, whose models are the divisible non-trivial ordered abelian groups.
\end{remark}

\begin{fact}[\cite{Elt.-prop.-of-ordered-abelian-groups,generalized_archimedean_groups,orderedfieldsdenseintheirrealclosure}]
For an ordered abelian group $G$ the following are equivalent:
\begin{enumerate}
\item $G$ is regular.
\item $G$ is elementarily equivalent to some archimedean ordered abelian group.
\item For any prime number $p$ and any infinite convex subset $A\subseteq G$ there is $b\in G$ such that $p\cdot b\in A$.
\end{enumerate}
If we assume in addition that $G$ is dense, then the above statements are moreover equivalent to the following:
\begin{enumerate}[resume]
\item $G$ is order-dense in its divisible hull.
\end{enumerate}
\end{fact}

\begin{fact}[\cite{Elt.-prop.-of-ordered-abelian-groups}]\label{convex Sg. elt. equiv.}\begin{enumerate}
    \item All regular discrete ordered abelian groups are elementarily equivalent to $(\mathbb{Z},+,<)$. 
    \item  Two non-trivial regular dense ordered abelian groups $G=(G,+,<)$ and $H=(H,+,<)$ are elementarily equivalent if and only if for any prime number $p$ we have that both $|G/pG|$ and $|H/pH|$ are infinite or $|G/pG|=|H/pH|$ holds. In particular a regular ordered abelian group is elementarily equivalent to any of its non-trivial convex subgroups.
\end{enumerate}

\end{fact}

As stated in \cite{defbilityHong} the following result follows for the equicharacteristic $0$ case already from a result of Koenigsmann in \cite{koenigsmann}.

\begin{fact}[\mbox{\cite[Theorem 4]{defbilityHong}}]\label{defbility_Hong}For any prime number $p$ there is a parameter-free $\mathcal{L}_{\mathrm{ring}}$-formula $\psi_{p}(x)$ such that for any field $K$ and any henselian valuation $v$ on $K$ with regular dense value group $\Gamma_v$ which is not $p$-divisible, one has $\psi_p(K)=\Ov$, i.e., the valuation $v$ is  defined by $\psi_p(x)$.
\end{fact}

We will use several times the following fact from general valuation theory.

\begin{fact}[\mbox{\cite[Theorem~4.4.2]{engler2005valued}}]\label{fact:comparable}
Let $v,w$ be henselian valuations with valuation rings $\mathcal{O}_{v}$, $\mathcal{O}_{w}$ on a field $K$, such that not both $k_v$ and $k_w$ are separably closed. Then $\Ov$ and $\Ow$ are comparable, i.e., $\Ov\subseteq\Ow$ or $\Ow\subseteq\Ov$ holds. 
\end{fact}

We assume the following result to be well-known but since we could not find it in the literature we will give a proof.

\begin{proposition}\label{uniqueness_of_valuation}
Let $v,w$ be henselian valuations with valuation rings $\mathcal{O}_{v}$, $\mathcal{O}_{w}$ on a field $K$, both with non-divisible regular dense value groups $\Gamma_{v}$, $\Gamma_{w}$ and residue fields $k_{v}$, $k_{w}$ of characteristic $0$. Then $\mathcal{O}_{v}=\mathcal{O}_{w}$ already holds.
\end{proposition}

\begin{proof}
First assume that $k_{v}$ or $k_{w}$ is not separably closed. Then $\Ov$ and $\Ow$ are comparable by Fact~\ref{fact:comparable}, i.e., $\Ov\subseteq\Ow$ or $\Ow\subseteq\Ov$. 
Let us assume w.l.o.g. that $\Ow\subseteq\Ov$. Then by \cite[Lemma~2.3.1]{engler2005valued} we have $w=\pi\circ v$, where $\pi:\Gamma_v\twoheadrightarrow\sfrac{\Gamma_{v}}{C}=\Gamma_{w}$ for some convex subgroup $C\subseteq \Gamma_{v}$. If $C$ were non-trivial, $\Gamma_{w}$ would be divisible by regularity of $\Gamma_{v}$, contradicting the assumptions.

Now assume $k_{v}$ and $k_{w}$ are both separably closed, thus algebraically closed since they are of characteristic $0$ by assumption. If $\Gamma_{v}$ and $\Gamma_{w}$ are not $p$-divisible for a common prime $p$, then $\Ov$ and $\Ow$ are defined by the same parameter-free formula, by Theorem~\ref{defbility_Hong}, whence $\Ov=\Ow$.
Otherwise let $p$ be a prime such that $\Gamma_{v}$ is $p$-divisible but $\Gamma_{w}$ is not. From the latter we infer that not every element of $K$ is a $p$-th power. On the other hand, using the former and the fact that $k_{v}$ is algebraically closed, an easy application of Hensel's Lemma yields that  every element of $K$ is a $p$-th power. This contradiction completes the proof.
\end{proof}

\begin{corollary}\label{elt.equiv&defb.val.}
Let $K,F$ be fields such that $K\equiv F$ in $\lring$. If there is a henselian valuation $v$ on $K$ with non-divisible regular dense value group $\Gamma_{v}$ and residue field $k_{v}$ of characteristic $0$, then there is exactly one henselian valuation $w$ on $F$ with non-divisible regular dense value group $\Gamma_{w}$ and residue field $k_{w}$ of characteristic $0$ and we have $\Gamma_{v}\equiv\Gamma_{w}$ and $k_{v}\equiv k_{w}$.
\end{corollary}

\begin{proof}
Fact~\ref{defbility_Hong} yields that the valuation $v$ is parameter-free definable whence $F$ has a valuation $w$ similarly defined. The uniqueness directly follows from Proposition~\ref{uniqueness_of_valuation}. Furthermore, by the interpretability of the residue field and of the value group, we get $k_{v}\equiv k_{w}$ and $\Gamma_v\equiv\Gamma_w$. 
\end{proof}

We proceed with some further observations that will prove useful later on.

\begin{notation}
Given a field $k$ and an ordered abelian group $\Gamma$ we denote by $k((t^\Gamma))$ the Hahn-series field with respect to $k$ and $\Gamma$. If not explicitly stated otherwise $\equiv$ will denote elementary equivalence in $\lring$ for fields and in $\lordgr$ for ordered groups. If we want to make this explicit, we will sometimes write $\equiv_{\mathrm{ring}}$ or $\equiv_{\mathrm{og}}$.
\end{notation}

Note that, when working with Hahn-series fields, one usually uses the additive notation for valuations. Since we use the multiplicative notation, we get, e.g., $t^1=1\in k((t^\Gamma))$. This abuse of notation should not lead to any confusion.

\begin{lemma}\label{div. Extension k((T))}
Let $\Gamma$ be a non-trivial regular ordered abelian group (this includes the divisible case) and $k$ a field of characteristic $0$. Then for any field $K$ which is elementarily equivalent to $k((t^\Gamma))$ in $\lring$, the following holds in $\mathcal{L}_{\mathrm{ring}}$:
\[k((t^\Gamma))\equiv K ((t^{\mathbb{R}^+})).
\]
\end{lemma}

\begin{proof}
By the AKE principle, we may assume that $K=k((t^\Gamma))$. As $ k((t^\Gamma)) ((t^{\mathbb{R}^+}))=k((t^{\Gamma\times\mathbb{R}^+}))$ where $\Gamma\times\mathbb{R}^+$ is equipped with the anti-lexicographical order, we naturally have a henselian valuation $v$ on $ k((t^\Gamma)) ((t^{\mathbb{R}^+}))$ with value group $\Gamma\times\mathbb{R}^+$. But then $(\Gamma\times\mathbb{R}^+)/\Gamma \cong \mathbb{R}^+
$ is divisible and $\Gamma$ is regular by definition and convex in the regular group $\Gamma\times\mathbb{R}^+$ so by Fact~\ref{convex Sg. elt. equiv.} it follows that $\Gamma\times\mathbb{R}^+\equiv \Gamma$. We apply the AKE principle and obtain
\[
\left(k((t^\Gamma)),\Gamma,k  \right)\equiv \left(k((t^\Gamma))((t^{\mathbb{R}^+})),\Gamma\times \mathbb{R}^+,k\right) 
\]
in the language of valued fields $\mathcal{L}_{\mathrm{val}}$, so in particular in $\mathcal{L}_{\mathrm{ring}}$.
\end{proof}

\begin{lemma}\label{simplification_lemma}
Let $k$ and $k^{\prime}$ be fields of characteristic $0$, and let $\Gamma$ and $\Gamma^{\prime}$ be non-trivial regular dense ordered abelian groups. If $k((t^\Gamma))\equiv_{\mathrm{ring}} k^{\prime}((t^{\Gamma^{\prime}}))$, then
\begin{enumerate}
    \item[\textup{(i)}] $\Gamma\equiv\Gamma^{\prime}$ and $k\equiv k^{\prime}$, or

\item[\textup{(ii)}] $k^{\prime}\equiv k((t^\Gamma))$ and $\Gamma^{\prime}$ is divisible, or 
\item[\textup{(iii)}] $k\equiv k^{\prime}((t^{\Gamma^{\prime}}))$ and $\Gamma$ is divisible.
\end{enumerate}
\end{lemma}

\begin{proof}
We may find, e.g., by the Keisler-Shelah Theorem, a common field extension $\Tilde{K}$ of $k((t^\Gamma))$ and $ k^{\prime}((t^{\Gamma^{\prime}}))$  and henselian valuations $v$ and $v^{\prime}$ on $\Tilde{K}$ such that $(\Tilde{K},v)\succcurlyeq_{\mathcal{L}_{\mathrm{val}}}k((t^{\Gamma}))$ and $(\Tilde{K},v')\succcurlyeq_{\mathcal{L}_{\mathrm{val}}}k'((t^{\Gamma'}))$. Note that $\Gamma\equiv\Gamma_v$ and $\Gamma'\equiv\Gamma_{v'}$, analogously for the residue fields involved.

\emph{Case 1: $\Gamma$ and $\Gamma'$ are both non-divisible.} Then $\Ov=\mathcal{O}_{v'}$ by Proposition~\ref{uniqueness_of_valuation}, so (i) holds.

\emph{Case 2: Both $\Gamma$ and $\Gamma'$ are divisible.} In particular, it follows that $\Gamma\equiv\Gamma'$. If $\Tilde{K}$ is algebraically closed, so are $k_v$ and $k_{v'}$, and thus (i) holds. Otherwise, $k_v$ and $k_{v'}$ are not separably closed, so $v$ and $v^{\prime}$ are comparable by Fact~\ref{fact:comparable}. If $\Ov=\mathcal{O}_{v'}$, we are in case (i). If $\Ov\subsetneq\mathcal{O}_{v'}$, there is a convex subgroup $(1)\subsetneq C\subsetneq \Gamma_v$ such that $\Gamma_{v'}\cong\Gamma_v/C$, and we get $k'\equiv k_{v'}\equiv k_v((t^C))\equiv k((t^\Gamma))$, so (ii) holds. Similarly, one shows that (iii) holds if $\Ov\supsetneq\mathcal{O}_{v'}$.

\emph{Case 3: Exactly one of $\Gamma'$ and $\Gamma$ is divisible.} W.l.o.g. we may assume that $\Gamma'$ is divisible and $\Gamma$ is not. Then $\Tilde{K}$ is not algebraically closed, so $k'\equiv k_{v'}$ is not separably closed (by the AKE principle, as $k'$ is of characteristic 0), and thus $v$ and $v'$ are comparable by Fact~\ref{fact:comparable}. As in Case 2, we infer that $\Gamma_{v'}\cong\Gamma_v/C$ for some convex subgroup $(1)\subsetneq C\subsetneq \Gamma_v$, as $\Gamma_{v}\cong\Gamma_{v'}/C'$ is impossible in this case. So again (ii) holds.
\end{proof}

\subsection{Continuous Logic}
In the following we will very briefly recall some notions from continuous logic but in general the reader shall be referred to \cite{yaacov_berenstein_henson_usvyatsov_2008} which will also be the main source for the rest of this section.
\begin{notation}
A \textit{continuous logic language} consists of non-logical symbols for predicates, functions and constants where the first two are all equipped with a fixed arity and a modulus of uniform continuity. Technically speaking, the non-logical part of the language also has to contain a positive real number $D$ which denotes an upper bound on the diameter of structures in this language and for each predicate $P$ a closed interval $I_{P}$ where $P$ takes its values. Throughout this text we will always assume $D=1$ and $I_{P}=[0,1]$ for any $P$.\\
As logical symbols we have a symbol for the metric $d(x,y)$ treated in a similar way as the equality symbol $=$ in classsical model theory, two quantifiers $\sup$ and $\inf$, as well as a set of connective symbols for the continuous functions $u:[0,1]^{n}\rightarrow [0,1]$.
Note that in general we do not have to allow the whole set of continuous functions $u:[0,1]^{n} \rightarrow [0,1]$ since it will be enough to \textit{uniformly approximate} formulas. We can even restrict ourselves to using finitely many connectives. See  \cite[Chapter~6]{yaacov_berenstein_henson_usvyatsov_2008} for a discussion of this topic.
\end{notation}

\begin{notation}
Given a continuous logic language $\mathcal{L}$ we define an \emph{$\mathcal{L}$-prestructure} to be a pseudo-metric space $(M_{0},d_{0})$ of diameter $\leq 1$ with interpretations of the functions, predicates and constants from $\mathcal{L}$. Here, each constant shall be interpreted by an element of $M_{0}$. Furthermore, each function symbol $f$ with arity $n$ shall be interpreted by a uniformly continuous map $M_{0}^{n}\rightarrow M_{0}$ having the modulus of uniform continuity as  specified by $f$. Finally, each predicate symbol $P$ with arity $m$ shall be interpreted by a uniformly continuous function $M_{0}^{m}\rightarrow [0,1]$ with modulus of uniform continuity as specified by $P$.
The associated \emph{$\mathcal{L}$-structure} will then be the completion of the quotient metric space of the prestructure where the interpretations of predicates and functions are in such a way that they extend those on the quotient metric space and are continuous. Hence, those interpretations are uniquely determined.

In the following we will denote by $M,N$ the underlying metric space of $\mathcal{L}$-structures $\mathcal{M},\mathcal{N}$ and often we will even identify $M$ with $\mathcal{M}$ (or $N$ with $\mathcal{N}$), i.e., the structure and its underlying metric space might be used interchangeably.
\end{notation}

\begin{remark}
All those notions can be easily generalised to multi-sorted structures and respective languages. Also, in a similar fashion as in the classical context we can inductively define terms, formulas, etc. and many concepts generalise to this setting. Especially, we can 
note that the continuous setting indeed includes the classical framework using the metric defined by $d(x,y):=0$ if $x=y$ and $d(x,y):=1$ otherwise and restricting predicate values to the set $\{0,1\}$.
\end{remark}

However, some concepts as definability (of sets) might not generalise completely intuitively at first glance. Again, for a full treatment we refer to \cite[Chapter~9]{yaacov_berenstein_henson_usvyatsov_2008}.

\begin{definition}[\mbox{\cite[9.1 and 9.16]{yaacov_berenstein_henson_usvyatsov_2008}}]
Given a continuous logic $\mathcal{L}$-structure $\mathcal{M}$, a subset $A\subseteq M$ and a uniformly continuous function $P:M^{n}\rightarrow [0,1]$ we say that $P$ is \textit{definable} in $\mathcal{M}$ over $A$ if there is a sequence $(\phi_{n}(x))_{n\in\mathbb{N}}$ of $\mathcal{L}(A)$-formulas such that the interpretations $\phi_{n}^{\mathcal{M}}(x)$ converge to $P(x)$ uniformly on $M^{n}$. In this case we call $P$ a \textit{definable predicate} (over $A$).

A closed set $D\subseteq M^{n}$ is a \textit{definable set} in $\mathcal{M}$ over $A$ if the distance (predicate) $\text{dist}(x,D)$ is definable in $\mathcal{M}$ over $A$.
\end{definition}

Next, we will focus on the ultraproduct construction in the continuous logic setting, namely \textit{metric ultraproducts}.
Metric ultraproducts have been studied and proved useful in several applications outside of model theory as well, for example in the context of Banach spaces \cite{ultraproduitsBanach} or in metric geometry in the build-up to a proof of Gromov's famous theorem on groups of polynomial growth \cite{gromov}.
We will give a short introduction to ultraproducts in continuous logic stating some definitions and results from \cite[Chapter~5]{yaacov_berenstein_henson_usvyatsov_2008}. Later we will deal with metric ultraproducts and classical logic ultraproducts at the
 same time. Thus, we will distinguish between them and label them either by $me$ or by $cl$.

\begin{notation}
Let $\mathcal{D}$ be an ultrafilter on the set $\mathcal{I}$, $X$ a topological space and let $(x_{i})_{i\in\mathcal{I}}$ be a sequence in $X$. Recall that $x$ is called an \textit{ultralimit} of $(x_{i})_{i\in\mathcal{I}}$ with respect to $\mathcal{D}$, denoted by $\lim_{i\rightarrow\mathcal{D}}x_{i}=x$, if for every neighbourhood $\mathcal{U}$ of $x$ the set $\{i\in\mathcal{I}\;|\;x_{i}\in\mathcal{U}\}$ is in $\mathcal{D}$. 
\end{notation}

\begin{fact}[\mbox{\cite[Lemma~5.1]{yaacov_berenstein_henson_usvyatsov_2008}}]
A topological space $X$ is compact Hausdorff if and only if for every $\mathcal{D}$ an ultralimit as above exists and is unique.
Given a continuous function $f:X\rightarrow X^{'}$ between topological spaces $X$ and $X'$ we have that
$\lim_{i\rightarrow\mathcal{D}}x_{i}=x$ implies $ \lim_{i\rightarrow\mathcal{D}}f(x_{i})=f(x).$
\end{fact} 

\begin{definition}\label{metricspaceUP}
Let $\mathcal{D}$ be an ultrafilter on $\mathcal{I}$ and $(M_{i},d_{i})_{i\in\mathcal{I}}$ a family of metric spaces, all with diameter $\leq 1$. Then there is a pseudo-metric on the cartesian product $\prod_{i\in\mathcal{I}}M_{i}$  defined by
$d(x,y)=\lim_{i\rightarrow\mathcal{D}} d_{i}(x_{i},y_{i})$, where $x=(x_{i})_{i\in\mathcal{I}}$ and $y=(y_{i})_{i\in\mathcal{I}}$. This naturally induces a metric on the quotient space $\sfrac{\prod_{i\in\mathcal{I}}M_{i}}{\sim_{\mathcal{D}}}$ where $x\sim_{\mathcal{D}}y$ if and only if $d(x,y)=0$. This space denoted by $\prod_{\mathcal{D}}^{\mathrm{me}}M_{i}$ is called the \emph{$\mathcal{D}$-ultraproduct} of $(M_{i},d_{i})_{i\in\mathcal{I}}$ and the corresponding equivalence classes are denoted by $((x_{i})_{i\in\mathcal{I}})_{\mathcal{D}}$ for $(x_{i})_{i\in\mathcal{I}}\in\prod_{i\in\mathcal{I}}M_{i}$.
\end{definition}

\begin{fact}[\mbox{\cite[Proposition~5.3]{yaacov_berenstein_henson_usvyatsov_2008}}]\label{UP_complete}
A metric ultraproduct of uniformly bounded complete metric spaces (as above) is always complete.
\end{fact}

Now, bearing in mind that in a continuous logic structure we have functions, predicates and constants given as uniformly continuous functions, for a family $(\mathcal{M}_{i})_{i\in\mathcal{I}}$ of metric structures in a given language we can define $\prod_{\mathcal{D}}^{\mathrm{me}}\mathcal{M}_{i}$, the \emph{$\mathcal{D}$-ultraproduct} of $(\mathcal{M}_{i})_{i\in\mathcal{I}}$, as a structure over the same language with the underlying space given by the metric space ultraproduct of the underlying metric spaces of the $\mathcal{M}_{i}$ and with functions, predicates and constants formed in the following manner:

Let $((M_{i},d_{i})_{i\in\mathcal{I}})$ and $((M_{i}^{'},d_{i}^{'})_{i\in\mathcal{I}})$ be families of metric spaces of diameter $\leq 1$ and (for a fixed $n\geq 1$) $f_{i}:M_{i}^{n}\rightarrow M_{i}^{'}$ uniformly continuous functions for each $i\in\mathcal{I}$ such that all have the same modulus of uniform continuity. Then we can define a function
$f_{\mathcal{D}}:\prod_{\mathcal{D}}^{\mathrm{me}}M_{i}^{n}\rightarrow \prod_{\mathcal{D}}^{\mathrm{me}}M_{i}^{'}$
that is still uniformly continuous with the same modulus of uniform continuity via:
\[f_{\mathcal{D}}\left(((x_{i}^{1})_{i\in\mathcal{I}})_{\mathcal{D}},\dots ,((x_{i}^{n})_{i\in\mathcal{I}})_{\mathcal{D}}\right)=((f_{i}(x_{i}^{1},\dots ,x_{i}^{n}))_{i\in\mathcal{I}})_{\mathcal{D}}.\]

(Note that we use here that the $\mathcal{D}$-ultrapower of the real interval $[0,1]$ can be identified with $[0,1]$ itself.) We also use that the underlying metric space of $\prod^{\mathrm{me}}_{\mathcal{D}}\mathcal{M}_i$ is complete by Fact~\ref{UP_complete}.

In the special case that $\mathcal{M}_i=\mathcal{M}$ for all $i$, we write $(\mathcal{M})_{\mathcal{D}}^{\mathrm{me}}$ for $\prod_{\mathcal{D}}^{\mathrm{me}}\mathcal{M}_{i}$, and we call it the \emph{$\mathcal{D}$-ultrapower of $\mathcal{M}$}.

By induction on the complexity of formulas we obtain the following equivalent of \L o\'s's Theorem in the metric setting. 

\begin{fact}[\mbox{\cite[Theorem~5.4 and 5.5]{yaacov_berenstein_henson_usvyatsov_2008}}]
Let $(\mathcal{M}_{i})_{i \in \mathcal{I}}$ be a family of $\mathcal{L}$-structures, $\mathcal{D}$ an ultrafilter on $\mathcal{I}$ and let $\mathcal{M}=\prod_{\mathcal{D}}^{\mathrm{me}}\mathcal{M}_i$. Then for every $\mathcal{L}$-formula $\varphi(x_{1}, \ldots, x_{n})$ and elements  $a_{k}=\left(\left(a_{i}^{k}\right)_{i \in \mathcal{I}}\right)_{\mathcal{D}}$ from $\mathcal{M}$, for $k=1, \ldots, n,$ one has
$$
\varphi^{\mathcal{M}}\left(a_{1}, \ldots, a_{n}\right)=\lim _{i\rightarrow\mathcal{D}} \varphi^{\mathcal{M}_{i}}\left(a_{i}^{1}, \ldots, a_{i}^{n}\right).
$$
Moreover, the diagonal embedding $\Delta:\mathcal{M}\rightarrow(\mathcal{M})^{\mathrm{me}}_{\mathcal{D}}$ is elementary.
\end{fact}

There is also a clean connection between ultraproducts of structures of a given class and the axiomatisability of that class. In addition we state the extension of the Keisler-Shelah Theorem to continuous logic which we will use later.

\begin{fact}[\mbox{\cite[Theorem~5.7]{yaacov_berenstein_henson_usvyatsov_2008}}]
If $\mathcal{M}$ and $\mathcal{N}$ are metric structures and $\mathcal{M}\equiv\mathcal{N}$, then there exists an ultrafilter $\mathcal{D}$ such that $(\mathcal{M})^{\mathrm{me}}_{\mathcal{D}}$ is isomorphic to $(\mathcal{N})^{\mathrm{me}}_{\mathcal{D}}$.
\end{fact}

\begin{fact}[\mbox{\cite[Proposition 5.14]{yaacov_berenstein_henson_usvyatsov_2008}}]\label{axiomatisability_criterium}
Suppose that $\mathcal{C}$ is a class of metric structures for some fixed language. Then $\mathcal{C}$ is axiomatisable if and only if $\mathcal{C}$ is closed under isomorphisms, ultraproducts and ultraroots. (Here, if $\mathcal{N}$ is an ultrapower of some $\mathcal{L}$-structure $\mathcal{M}$, then we call $\mathcal{M}$ an \textit{ultraroot} of $\mathcal{N}$.)
\end{fact}

\subsection{Metric valued fields}
In the following we will recall the basic notions from \cite{YAACOV_2014} which will allow us to consider certain valued fields as continuous logic structures. Precisely, we will restrict ourselves to metric valued fields. All results and notions in this subsection can be attributed to \cite{YAACOV_2014}.
\begin{notation}
A \textit{pre-metric valued field} is a valued field where the value group is a subgroup of the mutliplicative group of the positive real numbers. Let $K$ be such a field carrying a valuation $|\cdot|:K\rightarrow \mathbb{R}_{\geq0}$ then this naturally induces a metric on $K$ by $d(x,y):=|x-y|$. We will sometimes also write a pre-metric valued field as $K=(K,\Gamma,k)$ together with an ordered group embedding $\alpha:\Gamma\hookrightarrow \rplus$.

A \emph{metric valued field} is a complete pre-metric valued field.
\end{notation}
The problem we face if we want to consider a metric valued field as a continuous logic structure is that it is in general unbounded. Moreover, the straight-forward approach of working in a multi-sorted language containing sorts for closed balls of increasing radii does not work as well (see \cite[Proposition~1.2]{YAACOV_2014}). We will now briefly present  Ben Yaacov's idea to overcome this problem which is to work in the projective line over a metric valued field instead of working in the field itself. To recover the addition and multiplication from the field one uses a purely relational language with predicates for homogeneous polynomials. For the details of this approach and its necessity we refer to \cite{YAACOV_2014}.

If $f$ and $g$ are real valued functions, we will often write $f\vee g$ for $\max(f,g)$, and $f\wedge g$ for $\min(f,g)$.

\begin{notation}
Given a (pre-)metric valued field $K$, the \textit{projective line} $\kpro$ is the quotient $\sfrac{K^{2}\backslash\{0\}}{K^{\times}}$, whose elents are denoted by $[x:y]:=\sfrac{(x,y)}{K^{\times}}$. For any class we can find a representative $(x,y)$ with $|x|\lor |y|=1$ and thus the elements of $\kpro$ can be written as the classes $[x:y]$ where $|x|\lor |y|=1$. Note that this is of course not a unique representation since we can always multiply (both $x$ and $y$) with elements of $\{z\in K\;|\;|z|=1\}$. In order to simplify computations we will often even assume that our representatives are of the form $[x:1]$ or $[1:y]$. In general, elements of $\kpro$ shall be denoted by bold letters. Furthermore we fix to write $\am=[\ao:\as]$.
\end{notation}

\begin{definition}
Let $\X=(X_{1},\dots,X_{n}),\;\Xs=(X_{1}^{*},\dots,X_{n}^{*})$ and define $\Zhx\subseteq\Zxx$ to be the ring of polynomials which are homogeneous in each pair $(X_{i},X_{i}^{*})$ separately, i.e., polynomials $\Pxx$ such that for every $1\leq i\leq n$  there exists $r_{i}\in\mathbb{N}$ such that for every monomial $\PSxx$ of $\Pxx$ one has $\deg_{X_{i}}\PSxx + \deg_{X_{i}^{*}}\PSxx=r_{i}$.

The \textit{homogenisation} $\Phx\in\Zhx$ of a polynomial $\PX\in\Zx$ is then given by $\Phx:=P(\X/\Xs)\Ps$ 
where $\Ps=(\Xs)^{\deg_{\X}P}$, $\X/\Xs:=(X_1/X_1^*,\ldots,X_n/X_n^*)$ and $\deg_{\X}P=(\deg_{X_{1}}P,\dots,\deg_{X_{n}}P)$.
\end{definition}

\begin{definition}[\mbox{See \cite[Definition~1.4]{YAACOV_2014}}]
The language $\lp$ shall consist of predicates $\predpnxb$ for every polynomial $P_{n}(\X)\in\Zx$. The arity of $\predpnxb$ is $n$ where $P_{n}(\X)\in\Zxn$. For any $\predpnxb$ the modulus of uniform continuitiy shall be given by the identity. Additionally $\lp$ shall contain a constant symbol $\infty$.
\end{definition}

\begin{definition}[\mbox{See \cite[1.5]{YAACOV_2014}}]\label{projective_structure}
Let $\kpro$ be the projective line over a (pre-)metric valued field $(K,|\cdot|)$. We define an $\lp$-(pre-)structure on $\kpro$ by setting $\infty:=[1,0]$ and $\predpab:=|P^{h}(\Bar{\ao},\Bar{\as})|$ as well as $d(\am,\bmodel):=\|\am-\bmodel\|=|\ao\bs-\as\bo|$.
We will write $\predxs$ for the formula $d(x,\infty)$ and $\predpsxb$ for $\Pi\predxis^{\deg_{X_{i}}\!P}$.
\end{definition}

Note that $\kpro$ is an $\lp$-structure if and only if $K$ is a metric valued field.

\begin{fact}\label{F:Reconstruct}
There is an $\lp$-theory $\mathrm{MVF}$ (an explicit set of axioms is given in \cite[Definition~1.6]{YAACOV_2014}) such that a given $\lp$-structure is a model of $\mathrm{MVF}$ if and only if it is (canonically) isomorphic to $\kpro$ for some metric valued field $K$.

Moreover, extensions of models of $\mathrm{MVF}$ naturally correspond to embeddings of metric valued fields.
\end{fact}

\begin{proof}
This is proved in  \cite[Theorem~1.8]{YAACOV_2014}. The moreover part is not explicitly stated in \cite[Theorem~1.8]{YAACOV_2014}, but it easily follows from the proof given there.
\end{proof}


\begin{definition}
The $\lp$-theory of projective lines over  metric valued fields of equicharacteristic $0$ shall be denoted by $\mvf$. Moreover, the theory $\mvfd$ shall consist of $\mvf$ together with axioms stating that the value group is dense.
\end{definition}

\begin{remark}
The above is indeed axiomatizable in the continuous context by axioms expressing that $|p|=1$ for all prime numbers $p\in\mathbb{N}$ and $\inf_{x}|\|x\|-q|=0$ for all $q\in\mathbb{Q}\cap [0,1]$.
\end{remark}

\section{Ultraproducts and Residue Shift}\label{S:ult}

Now we will turn to metric ultraproducts of models of $\mvfd$. The aim is to understand them by the relation to the classical logic ultraproducts of their underlying metric valued fields. The difference between both ultraproducts  originates firstly from the fact that in the metric ultraproduct two sequences that are almost everywhere different can still give rise to the same element if the distances converge to zero. Secondly the value group in the metric formalism is bounded in its size, since it stays embedded in $\rplus$. Thus, given a sequence $x:=(x_{i})_{i\in\mathcal{I}}$ with $|x_{i}|<1$ $\mathcal{D}$-almost everywhere but $\lim_{i\rightarrow\mathcal{D}}|x_{i}|= 1$, the element $x$ will give rise to a \textit{new} element in the residue field of the metric ultraproduct, whereas it does not in the classical setting. This phenomenon possibly changes the elementary theory of the residue field in the metric ultraproduct. But still this change will turn out to be relatively tame and will be controlled by what we will call the \textit{residue shift}.

\begin{notation}
Throughout the rest of this section let $\mathcal{I}$ denote some index set and $\mathcal{D}$ an ultrafilter on $\mathcal{I}$ and $x:=(x_{i})_{i\in\mathcal{I}}$, $y:=(y_{i})_{i\in\mathcal{I}}$. For the moment we do not impose any further conditions on the ultrafilter $\mathcal{D}$. However the only case of interest will be that of a  countably incomplete (i.e., not closed under countable intersections) ultrafilter as justified by Lemma~\ref{Ginf_when_trivial}.
We further fix a family $\left(K_{i}\right)_{i\in \mathcal{I}}$ of  metric valued fields with valuations $v_{i}$, value groups $\Gamma_{i}\subseteq\rplus$ and residue fields $k_{i}$. Moreover we denote the embedding $\Gamma_{i}\hookrightarrow\rplus$ by $\alpha_{i}$.
\end{notation}

\begin{definition}\label{def_metr_UP}
The metric ultraproduct $K^{\mathrm{me}}$ is the underlying metric valued field of the structure $\upkipmetr$ which is the metric ultraproduct of the structures $(K_{i}\mathbb{P}^{1})_{i\in\mathcal{I}}$. We denote its residue field by $k^{\mathrm{me}}$ and its value group by $\Gamma^{\mathrm{me}}\subseteq\rplus$ and the valuation either by $v^{\mathrm{me}}$ or simply by $|\cdot|$.
\end{definition}
To relate $K^{\mathrm{me}}$ to its classical logic counterpart we have to fix some notations in the classical setting.
\begin{definition}\label{def_inf_group}
Let $\Gamma=\upgicl$ be the classical logic $\mathcal{D}$-ultraproduct taken in the language of ordered groups. In a canonical way, $\Gamma$ is an ordered subgroup of the classical logic ultrapower $\left(\mathbb{R}^+\right)_{\mathcal{D}}^{\mathrm{cl}}$ of the ordered group $\mathbb{R}^+$, with $\mathbb{R}^+\leq\left(\mathbb{R}^+\right)_{\mathcal{D}}^{\mathrm{cl}}$ diagonally embedded. Let $\Rinf$ be the subgroup of infinitesimals, i.e., the largest convex subgroup $\Delta$ of  $\left(\mathbb{R}^+\right)_{\mathcal{D}}^{\mathrm{cl}}$ such that $\Delta\cap\mathbb{R}^+=\{1\}$. Let $\Rfin$ be the the subgroup of finite elements, which is given by the convex hull of $\mathbb{R}^+$ in $\left(\mathbb{R}^+\right)_{\mathcal{D}}^{\mathrm{cl}}$. Set $\ginf:=\Gamma\cap\Rinf$ and $\gfin:=\Gamma\cap\Rfin$. 
\end{definition}

Note that $\ginf\leq\gfin$ are convex subgroups of $\Gamma$.

\begin{definition}\label{classicUP}
Let $\mathcal{K}_{i}$ be the classical logic $\lval$-structure with the underlying field $K_{i}$ and $\mathcal{K}:=\upKcalcl$ the classical logic ultraproduct whose underlying field is then given by $K:=\upKicl$ with valuation $v$ and value group $\Gamma:=\upgicl$ and residue field $k:=\upkicl$. Let $\ginf,\gfin$ be defined as above. Moreover, let $\bar{v}:K\rightarrow\sfrac{\Gamma}{\gfin}$ be the coarsening of $v$ and $\Bar{K}$ the corresponding residue field with induced valuation $\vfin:\bar{K}\rightarrow\gfin$.
Now let $\vfinb:\Bar{K}\rightarrow\sfrac{\gfin}{\ginf}$ be the coarsening of $\vfin$ and $\vinf:\bar{k}\rightarrow\ginf$ the induced valuation on the residue field of $\vfinb$. Additionally we set $\Bar{\Gamma}:=\gfininf$.
\end{definition}

\begin{lemma}\label{real_val_gr_classical}
There is a naturally induced embedding $\bar{\alpha}:\Bar{\Gamma}\hookrightarrow(\mathbb{R}^+,\cdot,1,<)$, i.e., we can assume $\Bar{\Gamma}\subseteq\mathbb{R}^{+}$.
\end{lemma}

\begin{proof}
Let $\alpha:\gfin\rightarrow\mathbb{R}^+$ be the standard part map. Then $\alpha$ is a group homomorphism with kernel $\ginf$. The induced map $\bar{\alpha}:\Bar{\Gamma}\rightarrow\mathbb{R}^+$ is easily seen to preserve $<$.
\end{proof}

\begin{theorem}\label{metr_UP_in_classical_UP}
The metric valued field $(K^{\mathrm{me}},v^{\mathrm{me}})$ is given by $(\Bar{K},\vfinb)$, in the sense that there is an isomorphism of valued fields $f:(\Bar{K},\vfinb)\rightarrow(K^{\mathrm{me}},v^{\mathrm{me}})$ that moreover induces the identity on $\bar\alpha(\bar{\Gamma})$.
\end{theorem}

\begin{proof}
We first show that there is a field isomorphism $g:\Bar{K}\rightarrow K^{\mathrm{me}}$. Given the quotient map $\Tilde{\beta}:\prodkip\rightarrow K^{\mathrm{me}}\mathbb{P}^{1}$, we have for $z=(z_{i})_{i\in\mathcal{I}}\in\prodkip$ that $\Tilde{\beta}(z)=\infty$ if and only if $\limiD\preds{z_{i}}=0$. Consequently, we obtain a map $\beta:\prodki\rightarrow K^{\mathrm{me}}\cup\{\infty\}$ with $\beta(x)=\infty$ if and only if $\gamma((v_i(x_i))_{i\in\mathcal{I}})>\gfin$ where $\gamma$ denotes the projection on the equivalence class given by the classical logic ultraproduct.  On the other hand consider the sequence
\[\prodki\xrightarrow{\gamma}K\xrightarrow{\plvb}\Bar{K}\cup\{\infty\}
\]where $\plvb$ is the place corresponding to $\bar{v}$. It follows from the definitions that  $\plvb\circ\gamma$ and $\beta$ are homomorphisms on the subring $\prodki\backslash Z$ of $\prodki$, where
$Z:=\left(\plvb\circ\gamma\right)^{-1}(\infty)=\beta^{-1}(\infty)$. Moreover we have that $\plvb\circ\gamma(x)=\plvb\circ\gamma(y)$ if and only if $\bar{v}(\gamma(xy^{-1}))<1\;\text{ in }\sfrac{\Gamma}{\gfin}$ (for $\gamma(x),\gamma(y)\in\mathcal{O}_{\bar{v}}$). Now the latter is equivalent to $\beta(x)=\beta(y)$ and it follows that $\Bar{K}\cong K^{\mathrm{me}}$ as fields.\\
It remains to show that $\vfinb$ and $v^{\mathrm{me}}$ define the same valuation on $\Bar{K}\cong K^{\mathrm{me}}$.  This directly follows from the definitions: Given $x\in\prodki$ we have
\[v^{\mathrm{me}}(\beta(x))=\limiD|x_{i}|=\limiD \alpha_{i}(v_{i}(x_{i}))=\alpha(v_{\mathrm{fin}}(\plvb\circ\gamma(x)))=\bar{\alpha}(\vfinb(\plvb\circ\gamma(x))).
\]
\end{proof}

\begin{lemma}\label{Ginf_when_trivial}
In the above setting, the following holds:
\begin{enumerate}
    \item If $\mathcal{D}$ is countably complete, then $\ginf$ is trivial and $\gfin=\Gamma$, and therefore $(K^{\mathrm{me}},\Gamma^{\mathrm{me}},k^{\mathrm{me}})=(K,\Gamma,k)$, canonically.
    \item If $\mathcal{D}$ is countably incomplete and the value groups $\Gamma_i$ are dense (and non-trivial) almost everywhere, then $\{1\}\subsetneq\ginf\subsetneq\gfin\subsetneq\Gamma$ and  $\Gamma^{\mathrm{me}}=\rplus$.
\end{enumerate}
\end{lemma}

\begin{proof}
We clearly have that $\ginf\supsetneq\{1\}$ if and only if there is a sequence $(x_{i})_{i\in\mathcal{I}}$ for $x_{i}\in K_{i}$ such that $\limiD|x_{i}|=1$ but $|x_{i}|< 1$ almost everywhere. 

Let us first prove (2), so we assume that $\mathcal{D}$ is countably incomplete and $\Gamma_i$ is dense (and non-trivial) almost everywhere. Then such a sequence exists (see, e.g., \cite[Lemma 10.59]{Drutu2018GeometricGT}). Moreover, given any $r\in\mathbb{R}^+$, as $\alpha_i(\Gamma_i)$ is dense in $\mathbb{R}^+$ for almost all $i$, for every $n>0$ we find $\gamma_n\in\Gamma$ such that $| r-\gamma_n|\leq1/n$. By $\aleph_1$-saturation of $\Gamma$, there is $\gamma\in \Gamma$ with $\gamma$ infinitesimally close to $r$. Then $\gamma\in\gfin$ and $\overline{\alpha}(\gamma \mod\ginf)=r$. This shows that $\Gamma^{\mathrm{me}}=\mathbb{R}^+$. Saturation also yields $\gfin\subsetneq\Gamma$ in this case.

To prove (1), we assume that $\mathcal{D}$ is countably complete. Then $\limiD|x_{i}|=1$ implies that $\|x_{i}\|=1$ 
almost everywhere (see \cite[Lemma 10.61]{Drutu2018GeometricGT}), from which it follows by what we said at the beginning of the proof that $\ginf$ is trivial. One shows similarly that $\gfin=\Gamma$  in this case. The result now follows from Theorem~\ref{metr_UP_in_classical_UP}.
\end{proof}

From now on we assume the ultrafilter $\mathcal{D}$ to be countably incomplete if not stated otherwise.

\begin{diagram}\label{diagram_of_valuations}
The above proof shows that the metric ultraproduct is given by the bottom line of this commutative diagram where the places are labeled with their corresponding valuations. Moreover it follows from 
\L o\'s's Theorem and general valuation theory (see, e.g.,  \cite[Corollary~4.1.4]{engler2005valued}) that all valuations occurring in the diagram are henselian.\\
\adjustbox{scale=1.1,center}{
\begin{tikzcd}
K \arrow[dd, "\bar{v}:K\rightarrow\Gamma/\Gamma_{\mathrm{fin}}" description] \arrow[rrrr, "v:K\rightarrow\Gamma"]                                                                 &  &  &  & k\cup\{\infty\}                                                                          & k \arrow[l, hook] \\
                                                                                                                                                                         &  &  &  &                                                                                          &                   \\
\bar{K}\cup\{\infty\}                                                                                                                                                    &  &  &  & \bar{k} \arrow[uu, "v_{\mathrm{inf}}:\bar{k}\rightarrow\Gamma_{\mathrm{inf}}" description] \arrow[d, hook] &                   \\
\bar{K} \arrow[u, hook] \arrow[rrrruuu, "v_{\mathrm{fin}}:\bar{K}\rightarrow \Gamma_{\mathrm{fin}}" description] \arrow[rrrr, "\bar{v}_{\mathrm{fin}}:\bar{K}\rightarrow\Gamma_{\mathrm{fin}}/\Gamma_{\mathrm{inf}}"] &  &  &  & \bar{k}\cup\{\infty\}                                                                    &                  
\end{tikzcd}}
\end{diagram}


\begin{remark}\label{ex_inf_val}
There is an induced valuation on the residue field of $K^{\mathrm{me}}$, given by $\vinf:\bar{k}\rightarrow\ginf$. We will call it the \textit{infinitesimal valuation}. Though this valuation exists, it is not captured by the valuation of the metric valued field, in other words, it is not captured by its metric. While this phenomenon leads to a possible change of the elementary theory of the residue field in an ultraproduct we can use that the infinitesimal value group has the same elementary theory as the value group of $K$ itself. As this value group is moreover regular it allows us to \textit{control} the elementary theory of the residue field. 

\end{remark}

\begin{proposition}\textbf{(Residue shift).}\label{residue_shift}
Assume that $(K_{i})_{i\in\mathcal{I}}$ is a family of equicharacteristic $0$ metric valued fields, and that $\mathcal{D}$ is a countably incomplete ultrafilter on $\mathcal{I}$, such that $\Gamma_{i}\equiv_{\mathrm{og}}\Delta$ for almost all $i\in\mathcal{I}$ and some fixed dense $\Delta\subseteq\rplus$. Furthermore, assume that $l\equiv_{\mathrm{ring}} k=\upkicl$. Then, $K^{\mathrm{me}}$ has value group $\Gamma^{\mathrm{me}}=\rplus$ and residue field $k^{\mathrm{me}}\equiv_{\mathrm{ring}} l((t^\Delta))$.
\end{proposition}

\begin{proof}
By Lemma~\ref{Ginf_when_trivial}(2), we have $\Gamma^{\mathrm{me}}=\rplus$  As we have seen in Remark~\ref{ex_inf_val}, $k^{\mathrm{me}}$ carries an infinitesimal valuation $\vinf: k^{\mathrm{me}}\rightarrow\ginf$, with residue field $k$, that is henselian and non-trivial by Lemma~\ref{Ginf_when_trivial}. Moreover, $\ginf\subseteq\Gamma$ is a non-trivial convex subgroup and since $\Delta\equiv\Gamma$ by \L o\'s's Theorem, $\Gamma$ is regular, so by Proposition~\ref{convex Sg. elt. equiv.} we get $\ginf\equiv\Gamma\equiv\Delta$. Now since $\vinf$ is henselian we can apply the classical logic AKE-principle and obtain that
\[\left(k^{\mathrm{me}},\ginf,k\right)\equiv\left(l((t^\Delta)),\Delta, l\right)
\]as valued fields and thereby in particular $k^{\mathrm{me}}\equiv_{\mathrm{ring}} l((t^\Delta))$.
\end{proof}

\section{Main theorem}\label{S:main}
In this section, in Theorem~\ref{metricAKEnew}, we will state the approximate AKE principle which describes the elementary classes (in $\lp$) of equicharacteristic $0$ metric valued fields in terms of the $\mathcal{L}_{\mathrm{og}}$-, resp. $\lring$-theories of their value groups and residue fields. As the residue shift (see Proposition~\ref{residue_shift}) showed, those theories are not necessarily preserved under taking the metric ultraproduct. However, we will see that this is the only obstruction. In other words, the elementary theory of a metric valued field as above is determined \textit{up to} the residue shift by the theories of value group and residue field. To make precise what we mean by \textit{up to} we will define (in Definition~\ref{definitionclassfomadmissiblepair}) a class of metric valued fields $\mathcal{C}(\mathfrak{G},\mathfrak{F})$ where $\mathfrak{G}$ and $\mathfrak{F}$ are theories of ordered groups and fields respectively. To properly define this class we will use the notion of an \textit{admissible pair}, given in Definition \ref{definitionadmissiblepair}. Its purpose is simply to make sure that $(\mathfrak{G},\mathfrak{F})$ is not a pair already obtainable by the residue shift from a different such pair. In the following definition, we will first describe those pairs of structures, called \textit{generating pairs}, whose theories yield an admissible pair.

\begin{definition}\label{def:generating-pair}
We say that a pair $(\Delta,l)$ consisting of a field $l$ of characteristic $0$ and a regular dense (non-trivial) ordered abelian group $\Delta$ is a \textit{generating pair} if one of the following two conditions is satisfied:
\begin{enumerate}[label=(\roman*)]
    \item $\Delta$ is not divisible;
    \item $\Delta$ is divisible  and $l\not\equiv l^{\prime}((t^{\Delta^{\prime}}))$ for any such pairs $(\Delta^{\prime},l^{\prime})$ with $\Delta^{\prime}$ non-divisible or with $\Delta^{\prime}$ divisible and $l^{\prime}\not\equiv l$.
\end{enumerate}
\end{definition}

\begin{definition}\label{definitionadmissiblepair}
    Let $\mathfrak{G}$ be a complete theory of non-trivial regular dense ordered abelian groups, $\mathfrak{F}$ a complete theory of fields of characteristic $0$ and let $(\Delta,l)$ be such that $\mathfrak{F}=\mathrm{Th}_{\mathrm{ring}}(l)$ and $\mathfrak{G}=\mathrm{Th}_{\mathrm{og}}(\Delta)$. We define $\mathfrak{F}^{\mathfrak{G}}:=\mathrm{Th}_{\mathrm{ring}}(l((t^{\Delta})))$ and $(\mathfrak{G},\mathfrak{F})^{\mathrm{sh}}:=(\mathrm{DOAG}, \mathfrak{F}^{\mathfrak{G}})$.
    Finally, we call the pair $(\mathfrak{G},\mathfrak{F})$ \textit{admissible}, if $(\Delta,l)$ is a generating pair. (Note that this does not depend on the choice of $(\Delta,l)$.)
\end{definition}
    
We now collect some properties of the \textit{(residue) shift-operator} $(\cdot)^{\mathrm{sh}}$. 
\begin{lemma}\label{lemmashiftoperator}
    Let $\mathfrak{G}$ and $\bar{\mathfrak{G}}$ be complete theories of non-trivial regular dense ordered abelian groups, and let $\mathfrak{F}$ and $\bar{\mathfrak{F}}$ be complete theories of fields of characteristic $0$. Then the following properties hold.
    \begin{enumerate}
        \item Idempotence: $((\mathfrak{G},\mathfrak{F})^{\mathrm{sh}})^{\mathrm{sh}}=(\mathfrak{G},\mathfrak{F})^{\mathrm{sh}}$.
        \item If $(\mathfrak{G},\mathfrak{F})$ is not admissible, then $(\mathfrak{G},\mathfrak{F})=(\mathfrak{G},\mathfrak{F})^{\mathrm{sh}}$. 
        \item If $(\mathfrak{G},\mathfrak{F})$ is admissible and $(\mathfrak{G},\mathfrak{F})=(\bar{\mathfrak{G}},\bar{\mathfrak{F}})^{\mathrm{sh}}$, then $(\mathfrak{G},\mathfrak{F})=(\bar{\mathfrak{G}},\bar{\mathfrak{F}})=(\mathrm{DOAG}, \mathfrak{F})=(\mathrm{DOAG},\mathfrak{F}^{\mathfrak{G}})$.
        \item If $(\mathfrak{G},\mathfrak{F})$ and $(\bar{\mathfrak{G}},\bar{\mathfrak{F}})$ are admissible and $(\mathfrak{G},\mathfrak{F})^{\mathrm{sh}}=(\bar{\mathfrak{G}},\bar{\mathfrak{F}})^{\mathrm{sh}}$, then $(\mathfrak{G},\mathfrak{F})=(\bar{\mathfrak{G}},\bar{\mathfrak{F}})$.
    \end{enumerate}
\end{lemma}

\begin{proof}
    Proof of (1): Lemma \ref{div. Extension k((T))} yields that $(\mathfrak{F}^{\mathfrak{G}})^{\mathrm{DOAG}}=\mathfrak{F}^{\mathfrak{G}}$.\\
    Proof of (2): If $(\mathfrak{G},\mathfrak{F})$ is not admissible, then $(\mathfrak{G},\mathfrak{F})=(\bar{\mathfrak{G}},\bar{\mathfrak{F}})^{\mathrm{sh}}$ for some $(\bar{\mathfrak{G}},\bar{\mathfrak{F}})\neq(\mathfrak{G},\mathfrak{F})$. Then, we can conclude by (1).\\
    Proof of (3): From $(\mathfrak{G},\mathfrak{F})=(\bar{\mathfrak{G}},\bar{\mathfrak{F}})^{\mathrm{sh}}$ it follows, that $\mathfrak{G}=\mathrm{DOAG}$ and $\mathfrak{F}=\bar{\mathfrak{F}}^{\bar{\mathfrak{G}}}$. Now, since $(\mathfrak{G},\mathfrak{F})$ is admissible, we get $\bar{\mathfrak{G}}=\mathrm{DOAG}$ and $\mathfrak{F}=\bar{\mathfrak{F}}$.\\
    Proof of (4): Follows from Lemma \ref{simplification_lemma} applied to models of $\mathfrak{F}^{\mathfrak{G}}=\bar{\mathfrak{F}}^{\bar{\mathfrak{G}}}$: Case (i) of Lemma \ref{simplification_lemma} directly implies the statement. In case (ii) or (iii), w.l.o.g. we get $\bar{\mathfrak{F}}=\mathfrak{F}^{\mathfrak{G}}$ and $\bar{\mathfrak{G}}=\mathrm{DOAG}$. Admissibility of $(\bar{\mathfrak{G}},\bar{\mathfrak{F}})$ then implies that $\mathfrak{G}=\mathrm{DOAG}$ and $\bar{\mathfrak{F}}=\mathfrak{F}$. Thus $(\mathfrak{G},\mathfrak{F})=(\bar{\mathfrak{G}},\bar{\mathfrak{F}})$.
\end{proof}

\begin{definition}\label{definitionclassfomadmissiblepair}
    Given an admissible pair $(\mathfrak{G},\mathfrak{F})$ we define the class $\mathcal{C}(\mathfrak{G},\mathfrak{F})$ to consist of all metric valued fields $\mathcal{K}=(K,\Gamma_{K},k_{K})$ such that for $\mathfrak{G}_{K}=\mathrm{Th}_{\mathrm{og}}(\Gamma_{K})$ and $\mathfrak{F}_{K}=\mathrm{Th}_{\mathrm{ring}}(k_{K})$ we have  $(\mathfrak{G},\mathfrak{F})=(\mathfrak{G}_{K},\mathfrak{F}_{K})$ or  $(\mathfrak{G},\mathfrak{F})^{\mathrm{sh}}=(\mathfrak{G}_{K},\mathfrak{F}_{K})$. For a generating pair $(\Delta,l)$ we set $\mathcal{C}(\Delta,l):=\mathcal{C}(\mathrm{Th}_{\mathrm{og}}(\Delta),\mathrm{Th}_{\mathrm{ring}}(l))$.
\end{definition}

\begin{remark}
    The classes $\mathcal{C}(\mathfrak{G},\mathfrak{F})$ are well defined and non-empty, and every metric valued field $K$ of equicharacteristic $0$ with dense value group belongs to a uniquely determined class of the above form.
\end{remark}
\begin{proof}
    The classes have empty intersection by Lemma~\ref{lemmashiftoperator} (2) and (3). Further, the pair $(\mathrm{Th}(\Gamma_{K}),\mathrm{Th}(k_{K}))$ is admissible itself or $(\mathrm{Th}(\Gamma_{K}),\mathrm{Th}(k_{K}))=(\mathrm{DOAG},\mathfrak{F}^{\mathfrak{G}})$ for some admissible pair $(\mathfrak{G},\mathfrak{F})$ and thus belongs to some class of the above form.
\end{proof}
We now state our main theorem which is a metric version of the Ax-Kochen-Ershov Theorem (later called \textit{metric AKE}) as elementary equivalence of two metric valued fields is reduced to elementary equivalence of residue field and value group defining the respective classes.
\begin{theorem}[Theorem A]\label{metricAKEnew}
Let $K,F$ be metric valued fields of equicha\-rac\-teristic $0$
 with dense value groups. Then $\kpro$ and $\fpro$ are elementarily equivalent if and only if they are in the same class $\mathcal{C}(\mathfrak{G},\mathfrak{F})$ for some admissible pair $(\mathfrak{G},\mathfrak{F})$.
\end{theorem}
\begin{examples}
We will now give several examples of classes $\mathcal{C}(\mathfrak{G},\mathfrak{F})$ for admissible pairs $(\mathfrak{G},\mathfrak{F})$. Let us say that a metric valued field $\mathcal{K}=(K,\Gamma_{K},k_{K})$ (of equicharacteristic 0 and with dense value group) is \textit{properly shifted}, if $(\Gamma_{K},k_K)$ is not a generating pair.
Note that if $\mathfrak{G}\neq \mathrm{DOAG}$ then the class $\mathcal{C}(\mathfrak{G},\mathfrak{F})$ contains properly shifted metric valued fields for any choice of $\mathfrak{F}$.

However, if $\mathfrak{G}=\mathrm{DOAG}$, two different cases can occur. Either $\mathcal{C}(\mathfrak{G},\mathfrak{F})$  contains properly shifted structures or $(\mathfrak{G},\mathfrak{F})^{\mathrm{sh}}=(\mathfrak{G},\mathfrak{F})$. The latter holds precisely, when $\mathfrak{F}^{\mathrm{DOAG}}=\mathfrak{F}$. We will call those classes \textit{fixed-point classes}. Many classes that naturally arise (e.g., all classes with residue field a local field of characteristic 0) are indeed fixed-point classes but we will also give examples of classes $\mathcal{C}(\mathfrak{G},\mathfrak{F})$ with $\mathfrak{G}= \mathrm{DOAG}$ that are not fixed-point classes.
\begin{itemize}
    \item \emph{Some fixed-point classes.}  For the following complete theories $\mathfrak{F}$ of fields of characteristic 0, $(\mathrm{DOAG},\mathfrak{F})$ is an admissible pair such that $\mathcal{C}(\mathrm{DOAG},\mathfrak{F})$ is a fixed-point class:
\begin{enumerate}
    \item $\mathfrak{F}=\mathrm{ACF}_0$. The corresponding class is that of all algebraically closed metric non-trivially valued fields of equicharacteristic 0.
    \item $\mathfrak{F}=\mathrm{RCF}$. The corresponding class is that of all real closed metric non-trivially valued fields, with convex valuation ring.
    \item $\mathfrak{F}=\mathrm{Th}_{\mathrm{ring}}(l)$ where $l$ a finite extension of $\mathbb{Q}_p$ for some prime $p$. The corresponding class is that of all $p$-adically closed metric valued fields, which are elementarily equivalent to $l$ in $\lring$ and such that the metric valuation is non-trivial and a proper coarsening of the $p$-adic valuation.
    \item $\mathfrak{F}=\mathrm{Th}_{\mathrm{ring}}(l)$ where $l=k((t))$ and $k$ is an arbitrary field of characteristic 0. Letting $\phi(x)$ be an $\lring$-formula such that $\phi(l)=k[[t]]$, the corresponding class is that of all metric valued fields, which are elementarily equivalent to $l$ in $\lring$ and such that the metric valuation is non-trivial and a proper coarsening of the valuation defined by $\phi(x)$. 
\end{enumerate}
The axiomatizability and completeness of the classes in 1. and 2. follow from Theorem~A, but they were already obtained in \cite{YAACOV_2014}, where the corresponding theories are further investigated. In 4., in order to show that $(\mathbb{R}^+,k((t)))$ is a generating pair, one may argue as in Case 3 in the proof of Lemma~\ref{simplification_lemma}. We leave the details to the reader.
\item \emph{Some non-fixed-point classes.} For the following complete theories $\mathfrak{F}$ of fields of characteristic 0, $(\mathrm{DOAG},\mathfrak{F})$ is an admissible pair such that $\mathcal{C}(\mathrm{DOAG},\mathfrak{F})$ is not a fixed-point class:
\begin{enumerate}
    \item $\mathrm{Th}_{\mathrm{ring}}(l)$ where $l$ non-large (e.g., $l$ any number field). Indeed, then $l\not\equiv l((t^{\mathbb{R}^+}))$, as $l((t^{\mathbb{R}^+}))$ is large and being large is first-order axiomatizable in $\lring$.\footnote{Let us sketch an elementary proof of the fact that $\mathbb{Q}\not\equiv\mathbb{Q}((t^{\mathbb{R}^+}))$. One may prove by elementary arguments that there are no $a,b \in\mathbb{Q}\backslash\{0\}$ such that $1+a^{3}=b^{3}$. On the other hand, the polynomial $P(X):=X^{3}-(1+t^{3})$ has a solution in $\mathbb{Q}((t^{\mathbb{R}^+}))$, by Hensel's Lemma.} (We refer to \cite{Pop14} for results on large fields.)
    \item $\mathrm{Th}_{\mathrm{ring}}(l)$ where $l$ is PAC and non-algebraically closed (e.g., $l$ any pseudofinite field). Then $l$ is large, but $l\not\equiv l((t^{\mathbb{R}^+}))$. Indeed, it follows from \cite[Theorem~10.14]{FrJa86} that $l((t^{\mathbb{R}^+}))$ is not PAC, which yields the result since being PAC is first-order axiomatizable in $\lring$. 
\end{enumerate}
\end{itemize}
\end{examples}

Remarkably the statement of the main theorem reduces now to elementary equivalence of metric valued fields seen as classical logic structures in $\lring$.

\begin{corollary}[Theorem B]\label{Cor:1storder}
Let $K,F$ be metric valued fields of equicha\-rac\-teristic $0$
 with dense value groups. Then $\kpro\equiv\fpro$ if and only if $K\equiv_{\mathrm{ring}} F$.
\end{corollary}

\begin{proof}
We have to show that $K\equiv_{\mathrm{ring}}F$ if and only if they are in the same class $\mathcal{C}(\mathfrak{G},\mathfrak{F})$.
If $K, \;F$ are in the same $\mathcal{C}(\mathfrak{G},\mathfrak{F})$, then this is a direct consequence of Lemma~\ref{div. Extension k((T))}, taking into account the definition of $\mathcal{C}(\mathfrak{G},\mathfrak{F})$.\\
For the other direction let $K\equiv_{\mathrm{ring}} F$ and the metric value groups on $K$ and $F$ shall be denoted by $\Gamma_{K}$ and $\Gamma_{F}$, the residue fields by $k_{K}$ and $k_{F}$. As $K$ and $F$ are metric valued fields, the following holds in $\lring$:
\[k_{K}((t^{\Gamma_{K}}))\equiv K\equiv F\equiv k_{F}((t^{\Gamma_{F}})).
\]Now, we can invoke Lemma~\ref{simplification_lemma} and directly conclude that $K$ and $F$ are in the same class $\mathcal{C}(\mathfrak{G},\mathfrak{F})$.
\end{proof}

As we have already mentioned, the discrete case is considerably easier. Given a metric valued field $K$, let us define the \textit{discreteness gap} $\dgk$ of $K$ to be
$\dgk:=\sup_{x\in K,\;|x|<1}|x|.$ Then $K$ is trivially valued if and only if $\dgk=0$, $K$ is discretely valued if and only if $0<\dgk<1$, and $K$ is non-trivially valued with dense value group if and only if $\dgk=1$. 

\begin{proposition}
Let $K, F$ be metric valued fields of equicharacteristic $0$ with discrete valuation. Then $\kpro\equiv F\mathbb{P}^{1}$ if and only if $\dgk=\dgf$ and $k_{K}\equiv k_{F}$. 
\end{proposition}

\begin{proof}
We only sketch the argument and leave the details to the reader. The discreteness gap is determined by the theory of a metric valued field. Moreover, in discrete metric valued fields the residue field is interpretable as an $\lring$- structure. This proves "$\Rightarrow$". For "$\Leftarrow$", note that in the discrete case Theorem~\ref{metr_UP_in_classical_UP} works similarly (with trivial infinitesimal valuation and exact same value group in the ultraproduct). Then we can conclude for example by applying Lemma~\ref{isomorphism-transfer} on some isomorphic ultrapowers.
\end{proof}

\section{Proof of Theorem~A}\label{S:proof}

The first goal is to prove a transfer for elementary equivalence of valued fields as classical logic structures to metric valued fields.
\begin{lemma}\label{isomorphism-transfer}
Let $(K_{1}, \Gamma_{1}, k_{1})$  and $(K_{2}, \Gamma_{2}, k_{2})$ be complete valued fields, and let $\alpha_{j}: \Gamma_{j}\hookrightarrow \rplus$ be embeddings, for $j=1,2$. Let $\sigma: K_{1}\cong K_{2}$ be an isomorphism of valued fields such that the induced isomorphism $\Tilde{\sigma}: \Gamma_{1}\cong\Gamma_{2}$ satisfies $\alpha_{2}\circ \Tilde{\sigma}=\alpha_{1}$. Then $\sigma$ induces an isomorphism $\sigma^{\mathrm{me}}:K_1\mathbb{P}^{1}\cong K_2\mathbb{P}^{1}$ of metric structures, where the metrics are induced by the embeddings $\alpha_{1}$ and $\alpha_2$.
\end{lemma}

\begin{proof}
Clear.
\end{proof}

\begin{proposition}\label{elt_equiv_transfer_final}
If two metric valued fields of equicharacteristic 0 have (full) value group $\rplus$ and are elementarily equivalent as classical logic structures in $\mathcal{L}_{\mathrm{val}}$, then their projective lines are elementarily equivalent as metric structures.
\end{proposition}
\begin{proof}
Let $K_{1}:=(K_{1}, \Gamma_{1}, k_{1})$ with $\alpha_{1}: \Gamma_{1}\cong \rplus$  and $K_{2}:=(K_{2}, \Gamma_{2}, k_{2})$ with $\alpha_{2}:\Gamma_{2}\cong \rplus$ be metric valued fields of equicharacteristic 0 with full value group. Assume that $K_{1}$ and $K_{2}$ are elementarily equivalent as classical logic structures in $\lval$.

The idea is to use the Keisler-Shelah Theorem to construct an isomorphism between classical logic ultrapowers of $K_{1}$ and $K_{2}$ that allows us to apply Lemma~\ref{isomorphism-transfer}. To do so we will work in the language $\lcval$ and choose constants $c_j\in \Gamma_{j}$ for $j=1,2$ such that $0<\alpha_{1}(c_{1})=\alpha_{2}(c_2)< 1$. As $(\Gamma_1,\cdot,<,c_1)\cong (\Gamma_2,\cdot,<,c_2)$, in particular we have  $(\Gamma_1,c_1)\equiv(\Gamma_2,c_2)$. Thus, $K_1\equiv K_2$ in $\lcval$ by Theorem~\ref{AKE}.

Now by Keisler-Shelah we find an ultrafilter $\mathcal{D}$ on an index set $\mathcal{I}$ such that there is an isomorphism $\sigma^{\mathrm{cl}}$ between the classical logic ultrapowers, in the language $\lcval$, of $K_1$ and $K_2$. Let $K_{j}^{\mathrm{me}}$ denote the underlying valued field of the metric ultrapowers $(\prod K_{j}\mathbb{P}^{1})_{\mathcal{D}}^{\mathrm{me}}$ for $j=1,2$.

We want to show that $\sigma^{\mathrm{cl}}$ induces an $\lval$-isomorphism $\sigma:K_{1}^{\mathrm{me}}\cong K_{2}^{\mathrm{me}}$ fulfilling the conditions of Lemma~\ref{isomorphism-transfer}, thus inducing an isomorphism $\sigma^{\mathrm{me}}:K_{1}^{\mathrm{me}}\mathbb{P}^1\cong K_{2}^{\mathrm{me}}\mathbb{P}^1$.
Using the relation between the classical and metric ultrapowers established in Proposition~\ref{metr_UP_in_classical_UP} it suffices that for the isomorphism $\sigma_{\Gamma}:\Delta_1:=\left(\prod\Gamma_{1}\right)_{\mathcal{D}}^{\mathrm{cl}}\cong \left(\prod\Gamma_{2}\right)_{\mathcal{D}}^{\mathrm{cl}}=:\Delta_2$ induced by $\sigma^{\mathrm{cl}}$ we have that $\sigma_{\Gamma}(c_{1})=c_{2}$ and  $\sigma_{\Gamma}(\Delta_{1,\mathrm{fin}})=\Delta_{2,\mathrm{fin}}$. The former is clear, since $\sigma^{\mathrm{cl}}$ is an $\lcval$-isomorphism. The latter follows from the former, as, by definition of the constants $c_{j}$ and of $\Delta_{j,\mathrm{fin}}$, we have that $\Delta_{j,\mathrm{fin}}$ is in both cases given as the smallest convex subgroup of $\Delta_j$ that contains $c_j$. 
Consequently we obtain $K_{1}\mathbb{P}^{1}\equiv K_{2}\mathbb{P}^{1}$ as metric structures in $\lp$.
\end{proof}

We now restate the residue shift (Proposition \ref{residue_shift}) using the terminology introduced in the last chapter.

\begin{corollary}\textbf{(Residue shift).}\label{corollaryresidue_shift}
Assume that $(\mathfrak{G},\mathfrak{F})$ is a pair of complete theories of some regular dense group and some field of characteristic $0$. Further assume that $(K_{i})_{i\in\mathcal{I}}$ is a family of equicharacteristic $0$ metric valued fields, and that $\mathcal{D}$ is a countably incomplete ultrafilter on $\mathcal{I}$, such that $\Gamma_{i}\models\mathfrak{G}$ and $k_{i}\models \mathfrak{F}$ for almost all $i\in\mathcal{I}$. Then, for the value group and residue field of the metric ultraproduct $K^{\mathrm{me}}$ we have $\Gamma^{\mathrm{me}}=(\mathbb{R}^+,\cdot)$ and $k^{\mathrm{me}}\models \mathfrak{F}^{\mathfrak{G}}$, in particular, $(\mathrm{Th}_{\mathrm{og}}(\Gamma^{\mathrm{me}}),\mathrm{Th}_{\mathrm{ring}}(k^{\mathrm{me}}))=(\mathfrak{G},\mathfrak{F})^{\mathrm{sh}}$.
\end{corollary}

\begin{corollary}\label{elt_equiv_transfer_for_shifted}
If two metric valued fields $K_{1}$ and $K_{2}$ are in the same class associated to an admissible pair, then $K_{1}\mathbb{P}^1$ and $K_{2}\mathbb{P}^1$  are elementarily equivalent.
\end{corollary}

\begin{proof}
Let $K\in \mathcal{C}(\mathfrak{G},\mathfrak{F})$ where $(\mathfrak{G},\mathfrak{F})$ is an admissible pair. Then either $(\mathfrak{G}_K,\mathfrak{F}_K)=(\mathfrak{G},\mathfrak{F})$ or $(\mathfrak{G}_K,\mathfrak{F}_K)=(\mathfrak{G},\mathfrak{F})^{\mathrm{sh}}$. Take a countably incomplete ultrafilter $\mathcal{D}$ on an index set $\mathcal{I}$. Then by Corollary~\ref{corollaryresidue_shift} (applied to $(\mathfrak{G}_K,\mathfrak{F}_K)$) and using, if necessary, that $(\mathfrak{G},\mathfrak{F})^{\mathrm{sh}}=((\mathfrak{G},\mathfrak{F})^{\mathrm{sh}})^{\mathrm{sh}}$ we get that $K^{\mathrm{me}}$ has full value group $(\mathbb{R}^+,\cdot)$ and 
 \[(\mathfrak{G}_{K^{\mathrm{me}}},\mathfrak{F}_{K^{\mathrm{me}}})=(\mathfrak{G}_K,\mathfrak{F}_K)^{\mathrm{sh}}=(\mathfrak{G},\mathfrak{F})^{\mathrm{sh}}.\]
Therefore, by the classical AKE-principle (Theorem~\ref{AKE}), for $K_1,K_2\in \mathcal{C}(\mathfrak{G},\mathfrak{F})$ the valued fields $K_{1}^{\mathrm{me}}$ and $K_{2}^{\mathrm{me}}$ are elementarily equivalent as structures in $\mathcal{L}_{\mathrm{val}}$ and as both have full value group we can apply Proposition~\ref{elt_equiv_transfer_final} to obtain $K_{1}^{\mathrm{me}}\mathbb{P}^1\equiv K_{2}^{\mathrm{me}}\mathbb{P}^1$. Now using \L o\'s's Theorem in $\lp$ we have $K_{1}\mathbb{P}^1\equiv K_{1}^{\mathrm{me}}\mathbb{P}^1$ and $K_{2}\mathbb{P}^1\equiv K_{2}^{\mathrm{me}}\mathbb{P}^1$ and so  $K_{1}\mathbb{P}^1\equiv K_{1}^{\mathrm{me}}\mathbb{P}^1\equiv K_{2}^{\mathrm{me}}\mathbb{P}^1\equiv K_{2}\mathbb{P}^1,$
which completes the proof.
\end{proof}

To finalise the proof of Theorem~\ref{metricAKEnew} it now only remains to show that the classes are indeed elementary. We will invoke Proposition~\ref{axiomatisability_criterium}, hence it suffices to show that the classes are closed under taking ultraproducts and ultraroots.

\begin{lemma}
The classes defined in the main theorem are closed under taking ultraproducts.
\end{lemma}

\begin{proof}
Let $(K_{i})_{i\in\mathcal{I}}$ be a family of  metric valued fields from the same class $\mathcal{C}(\mathfrak{G},\mathfrak{F})$, and let $\mathcal{D}$ be an ultrafilter on $\mathcal{I}$. 
As before we denote the underlying valued field of the metric ultraproduct by $(K^{\mathrm{me}},\Gamma^{\mathrm{me}},k^{\mathrm{me}})$, and the classical logic ultraproduct of the family $(K_{i})_{i\in\mathcal{I}}$ by $(K,\Gamma,k)$. If $\mathcal{D}$ is countably complete, 
$(K^{\mathrm{me}},\Gamma^{\mathrm{me}},k^{\mathrm{me}})=(K,\Gamma,k)$ by Lemma~\ref{Ginf_when_trivial}(1), so $K^{\mathrm{me}}\in\mathcal{C}(\mathfrak{G},\mathfrak{F})$ by \L o\'s's Theorem. If $\mathcal{D}$ is countably incomplete we can apply Corollary \ref{corollaryresidue_shift} as in the proof of Corollary~\ref{elt_equiv_transfer_for_shifted}.\end{proof}

\begin{lemma}\label{going_down}
The classes defined in the main theorem are closed under taking ultraroots.
\end{lemma}

\begin{proof}
Let $K_0=(K_0,\Gamma_0,k_0)$ be a metric valued field and $\mathcal{D}$ an ultrafilter on some set. Let $K=(K,\Gamma,k)=(K_0)^{\mathrm{cl}}_{\mathcal{D}}$ in $\lval$ and let $(K^{\mathrm{me}},\Gamma^{\mathrm{me}},k^{\mathrm{me}})$ be the metric valued field such that $K^{\mathrm{me}}\mathbb{P}^1=(K_0\mathbb{P}^1)^{\mathrm{me}}_{\mathcal{D}}$ in $\lp$. Assume that $K^{\mathrm{me}}\in \mathcal{C}(\mathfrak{G},\mathfrak{F})$ for some admissible pair $(\mathfrak{G},\mathfrak{F})$. We need to show that $K_0\in\mathcal{C}(\mathfrak{G},\mathfrak{F})$. If $\mathcal{D}$ is countably complete, then $(K^{\mathrm{me}},\Gamma^{\mathrm{me}},k^{\mathrm{me}})=(K,\Gamma,k)$  by Lemma~\ref{Ginf_when_trivial}(1). As $K_0\equiv K$ in $\lval$, it then  follows that $K_0\in\mathcal{C}(\mathfrak{G},\mathfrak{F})$. From now on $\mathcal{D}$ is assumed to be countably incomplete.\\
By Corollary \ref{corollaryresidue_shift} we have $(\mathrm{Th}_{\mathrm{og}}(\Gamma^{\mathrm{me}}),\mathrm{Th}_{\mathrm{ring}}(k^{\mathrm{me}}))=(\mathrm{Th}_{\mathrm{og}}(\Gamma_{0}),\mathrm{Th}_{\mathrm{ring}}(k_{0}))^{\mathrm{sh}}$. If $(\mathrm{Th}_{\mathrm{og}}(\Gamma^{\mathrm{me}}),\mathrm{Th}_{\mathrm{ring}}(k^{\mathrm{me}}))=(\mathfrak{G},\mathfrak{F})$, then $(\mathrm{Th}_{\mathrm{og}}(\Gamma^{\mathrm{me}}),\mathrm{Th}_{\mathrm{ring}}(k^{\mathrm{me}}))$ is admissible and by Lemma \ref{lemmashiftoperator} (2) it follows that $(\mathrm{Th}_{\mathrm{og}}(\Gamma_{0}),\mathrm{Th}_{\mathrm{ring}}(k_{0}))=(\mathfrak{G},\mathfrak{F})$.\\ 
If $(\mathrm{Th}_{\mathrm{og}}(\Gamma^{\mathrm{me}}),\mathrm{Th}_{\mathrm{ring}}(k^{\mathrm{me}}))=(\mathfrak{G},\mathfrak{F})^{\mathrm{sh}}$, we have $(\mathrm{Th}_{\mathrm{og}}(\Gamma_{0}),\mathrm{Th}_{\mathrm{ring}}(k_{0}))^{\mathrm{sh}}=(\mathfrak{G},\mathfrak{F})^{\mathrm{sh}}$. Now, if $(\mathrm{Th}_{\mathrm{og}}(\Gamma_{0}),\mathrm{Th}_{\mathrm{ring}}(k_{0}))$ is admissible, then we can conclude by Lemma \ref{lemmashiftoperator} (3). Otherwise, if $(\mathrm{Th}_{\mathrm{og}}(\Gamma_{0}),\mathrm{Th}_{\mathrm{ring}}(k_{0}))$ is not admissible, by Lemma \ref{lemmashiftoperator} (2) we have $(\mathrm{Th}_{\mathrm{og}}(\Gamma_{0}),\mathrm{Th}_{\mathrm{ring}}(k_{0}))^{\mathrm{sh}}=(\mathrm{Th}_{\mathrm{og}}(\Gamma_{0}),\mathrm{Th}_{\mathrm{ring}}(k_{0}))$ and thus it follows that $(\mathrm{Th}_{\mathrm{og}}(\Gamma_{0}),\mathrm{Th}_{\mathrm{ring}}(k_{0}))=(\mathfrak{G},\mathfrak{F})^{\mathrm{sh}}$. In both cases, we have $K_{0}\in\mathcal{C}(\mathfrak{G},\mathfrak{F})$.  
\end{proof}

\section{Metric valued difference fields}\label{S:diff}
In this last section, we will prove Theorem C on the non-existence of a model-companion for the theory of metric valued difference fields, thus answering a question of Ben Yaacov negatively. 

Before we get to the proof, we will put our work into a larger context, recalling some results on isometric valued difference fields in the classical context, where a model-companion does exist in equicharacteristic 0.

\subsection{Isometric valued difference fields in the classical context}
An \emph{isometric valued difference field} is a valued field $(K,\Gamma_K,k_K)$ together with an isomorphism $\sigma$ that induces the identity on the value group $\Gamma_K$. We denote the induced $\lring$-automorphism of $k_K$ by $\overline{\sigma}$, and we consider isometric valued difference fields in the language $\lvals$ given by $\lval$ augmented by function symbols for $\sigma$ and for $\overline{\sigma}$. 

Let $T_{\mathrm{val}}^{\mathrm{iso}}$ be the theory of (henselian) isometric valued difference fields in equicharacteristic 0, considered in the language $\lvals$.

\begin{fact}[\cite{frobeniuswitt}]\label{F:BMS}
The theory $T_{\mathrm{val}}^{\mathrm{iso}}$ admits a model-compantion $\mathrm{\mathrm{VFA}}^{\mathrm{iso}}$, which may be axiomatised as follows: For $\mathcal{K}=(K,\Gamma_K,k_K,\sigma)\models T_{\mathrm{val}}^{\mathrm{iso}}$, one has  $\mathcal{K}\models \mathrm{VFA}^{\mathrm{iso}}$ if and only if the following conditions hold:
\begin{enumerate}
    \item $\mathcal{K}$ has \emph{enough constants}, i.e.,  $\Gamma_K=\Gamma_{Fix(\sigma)}$.
    \item $\Gamma_K\models \mathrm{DOAG}$
    \item $(k_K,\overline{\sigma})\models \mathrm{ACFA}$
    \item $\mathcal{K}$ is $\sigma$-henselian\footnote{See \cite{frobeniuswitt} for the definition of $\sigma$-henslianity.}.
\end{enumerate}
\end{fact}

As both DOAG and ACFA are NTP$_2$ theories, it follows from \cite[Theorem~4.6]{ChHi14} that any completion of $\mathrm{VFA}^{\mathrm{iso}}$ is NTP$_2$. Moreover,  $\mathrm{VFA}^{\mathrm{iso}}$ is arithmetically meaningful. Indeed, for $p$ a prime number let $v_p$ be the $p$-adic valuation on $\mathbb{C}_p$, and let $\sigma_p$ be an isometric lift of the Frobenius automorphism on $k_{\mathbb{C}_p}=\mathbb{F}_p^{alg}$. 
Using Hrushovski's deep characterization of the non-standard Frobenius automorphism from \cite{Hru04}, one may infer in an elementary way from  Fact~\ref{F:BMS} (see \cite[Section~2.4]{HiRi24}) the following characterisation of the $\lvals$-theory $\mathrm{VFA}^{\mathrm{iso}}$.

\begin{fact}
$\mathrm{VFA}^{\mathrm{iso}}=\{\phi\mid (\mathbb{C}_p,v_p,\sigma_p)\models\phi\text{ for all $p\gg0$}\}$.
\end{fact}





\subsection{Metric valued difference fields}
We now get back to the continuous setting. Again we will work in the projective line rather than in the field itself. We expand the language $\lp$ to a language $\lps$ that consists of all of $\lp$ together with an additional function symbol $\sigma$ having the identity as modulus of uniform continuity (since $\sigma$ will be isometric).
\begin{definition}
The theory $\mathrm{MVF}_\sigma$ (metric valued fields with automorphism) shall consist of the following set of axioms:
\begin{align*}
 &(I)&  &\mathrm{MVF} \\
 &(II)&  &\|P(\bar{x})\|=\|P(\sigma(\bar{x}))\|&&\text{for any }\;P(\bar{X})\in\mathbb{Z}[\bar{X}] \\
 &(III)& &d(\infty,\sigma(\infty))=0 \\
 &(IV)& &\sup\limits_{y}\inf\limits_{x}\;\|\sigma(x)-y\|=0\\
 \end{align*}
\end{definition}

\begin{proposition}
\begin{enumerate}
    \item Any metric valued field $K$ endowed with an isometric automorphism $\Tilde{\sigma}$ gives rise to an $\lps$-structure $(K\mathbb{P}^1,\sigma)$ which is a model of $\mathrm{MVF}_\sigma$, by setting, for any $\am=[\ao:\as]$, $\sigma(\am):=[\Tilde{\sigma}(\ao):\Tilde{\sigma}(\as)]$. 
    \item Any model of $\mathrm{MVF}_\sigma$ arises in this way, and the models of $\mathrm{MVF}_\sigma$ are precisely the models of $\mathrm{MVF}$ endowed with an automorphism.
\end{enumerate}
\end{proposition}

\begin{proof}
1. is straight forward. 

To prove 2., it is enough to show that for any model $(K\mathbb{P}^1,\sigma)\models \mathrm{MVF}_\sigma$, the function $\sigma$ is an automorphism of the metric structure $K\mathbb{P}^1$, as it then automatically comes from an isometric automorphism of $K$, by the moreover part of Fact~\ref{F:Reconstruct}.

By $(II)$ and $(III)$, $\sigma$ is an isometric self-embedding of $K\mathbb{P}^1$, noting that $d(x,y)=\|P_d(x,y)\|$, where $P_d(x,y)=x-y$, whereby $(II)$ forces $\sigma$ to be an isometry. 

Now, $\sigma(K\mathbb{P}^1)$ is dense in $K\mathbb{P}^1$ by $(IV)$, and it is a complete  metric subspace. Hence $\sigma(K\mathbb{P}^1)=K\mathbb{P}^1$, and so $\sigma$ is surjective. This finishes the proof.
\end{proof}

\begin{definition}
We call a pair $(K,\Tilde{\sigma})$ as above a \textit{metric valued difference field}. Further we will reduce to using $\sigma$ for both the valued field automorphism on $K$ and the respective structure automorphism on $K\mathbb{P}^{1}$.
\end{definition}

The following lemma is the key ingredient in our proof of the non-existence of a model companion.

\begin{lemma}\label{L:phi}
Let $\mathcal{M}=(K\mathbb{P}^1,\sigma)\models \mathrm{MVF}_\sigma$ and $\am=[\ao:\as] \in K\mathbb{P}^1$. Consider the  $\lps$-formula \[\phi(x)=\inf\limits_{y}\left[ \|yx-\sigma(y)\|+\left((1-\|y\|)\vee (1-\|y^{*}\|)\right) \right].\]
Then the following holds:
\begin{enumerate}
    \item 
    If  $\n{\am}\neq 1$ or $\sn{\am}\neq 1$, then $\phi^{\mathcal{M}}(\am)=1$.
    \item  If $\am$ corresponds to an element $a\in K^\times$ such that $a=\sigma(b)/b$ for some $b\in K^{\times}$ with $| b|=1$ (so in particular $\n{\am}= \sn{\am}=1$), then $\phi^{\mathcal{M}}(\am)=0$.
\end{enumerate}
\end{lemma}

\begin{proof}
To prove 1. it suffices to show that for every $\bmodel\in K\mathbb{P}^1$ we have
\[\n{\bmodel\am-\sigma(\bmodel)}=\n{\bmodel}\wedge\sn{\bmodel}. \]
Since $\sigma$ is valuation-preserving, $|\first{\sigma(b)}|=|\bo|$ and $|\second{\sigma(b)}|=|\bs|$ hold, and as by assumption either $|\ao|<1$ or $|\as|<1$ we get,  by using the ultra-metric inequality, $\|\bmodel\am-\sigma(\bmodel)\|=|\bo\ao\sigma(\second{b})-\sigma(\bo)\bs\as|=|\bo\bs|=\|\bmodel\|\wedge\predbs$.

The proof of 2. is clear.
\end{proof}

The next aim will be to find for every model $\mathcal{M}\models \mathrm{MVF}_\sigma$ and every $\am\in\mathcal{M}$ such that $\n{\am}= \sn{\am}=1$ a model $\mathcal{M}\subseteq\mathcal{N}\models \mathrm{MVF}_\sigma$ such that $\phi(\am)=0$ in $\mathcal{N}$.

\begin{lemma}\label{extensionconstruction}
Given a metric valued difference field $(K,\sigma)$ and $\af\in K$ with $|\af |=1$ there exists a metric valued difference field $(F,\Tilde{\sigma})$ extending $(K,\sigma)$ and $\bfield\in F$ such that $|\bfield |=1$ and  $\af=\sigma(\bfield)/\bfield$.
\end{lemma}

\begin{proof}
We directly construct $(F,\Tilde{\sigma})$. Let $A=K\left[ X \right]$ and set $\left|\sum\limits_{i=0}^{n}\cf X^{i} \right|=\max\limits_{0\leq i\leq n}|\cf_{i}|$ for any such polynomial. Further let $\Tilde{\sigma}\restriction_{K}=\sigma$ and $\Tilde{\sigma}\left(\sum\limits_{i=0}^{n}\cf X^{i} \right)=\sum\limits_{i=0}^{n}\sigma(\cf)\af^{i} X^{i} $.

The metric valuation we obtain on $K(X)$ is the Gauss extension of $K$. Let $F$ be the completion of $(K(X),\left|\cdot \right|)$. Then $\Tilde{\sigma}$ extends to an isometric automorphism of $F$. Moreover, $\Tilde{\sigma}(X)/X=\af$ and $| X |=1$, which completes the proof. 
\end{proof}



\begin{theorem}[Theorem C]\label{thm:no-model-companion}
Fix any $(a,b)\in\{(0,0),(0,p),(p,p)\mid p\text{ prime}\}$. Then the theory $\mathrm{MVF}_{(a,b),\sigma}$ does not have a model-companion.
\end{theorem}

\begin{proof}
Assume otherwise and let $T$ denote the theory of the model companion. Since  $\mathrm{MVF}_{(a,b),\sigma}$ is an inductive theory, it admits existentially closed models, and the models of $T$ are precisely the existentially closed models of $\mathrm{MVF}_{(a,b),\sigma}$. 

Let now $\mathcal{M}=(K\mathbb{P}^{1},\sigma)$ be an existentially closed model of $\mathrm{MVF}_\sigma$. 

\begin{claim}\mbox{}
\begin{enumerate}
    \item The value group $\Gamma_K$ is a dense subgroup of $\mathbb{R}^+$.
    \item For any $\am\in K\mathbb{P}^1$ with $\sn{\am} =\n{\am}= 1$ one has $\phi^{\mathcal{M}}(\am)=0$.
\end{enumerate}
 
\end{claim}

Indeed, it easy to see that every metric valued difference field extends to an algebraically closed non-trivially valued one, so in particular to one with value group a dense subgroup of $\mathbb{R}^+$. 
Combining Lemma~\ref{extensionconstruction} with Lemma~\ref{L:phi}(2), we obtain the second assertion, thus proving the claim.

We now consider the partial type over $\emptyset$ given by \[\pi(x):=\{\phi(x)=1\}\cup\{(1-\|x\|)\vee (1-\|x^{*}\|)\leq1/n:n\in \mathbb{N}_{>0}\}.\]  Then $\pi$  is finitely satisfiable in $\mathcal{M}$ by Lemma~\ref{L:phi}(1) and the first part of the claim. Thus, in some elementary extension   $\mathcal{N}$ of $\mathcal{M}$ there exists $\am\models\pi$. As $\mathcal{N}\models T$, in particular $\mathcal{N}$ is existentially closed. But, in $\mathcal{N}$, we have $\sn{\am} =\n{\am}= 1$ and $\phi(\am)=1$,  contradicting the second part of the claim.
\end{proof}

\bibliography{MVF}
\bibliographystyle{abbrv}
\end{document}